\documentclass[reqno,11pt]{amsart} 
\usepackage{amsmath}
\usepackage{amssymb}
\usepackage{amsthm}

\newcommand\NoBlackBoxes{\global\overfullrule0pt}
\NoBlackBoxes
\theoremstyle{plain} % default

\def\4{\kern1pt}

\def\6{\vphantom0}

\def\8{\kern-10pt}
\def\7#1{_{(#1)}}

\begin{document}

\def\ffrac#1#2{\raise.5pt\hbox{\small$\4\displaystyle\frac{\,#1\,}{\,#2\,}\4$}}
\def\ovln#1{\,{\overline{\!#1}}}
\def\ve{\varepsilon}
\def\kar{\beta_r}

\title{NON-UNIFORM BOUNDS IN THE POISSON  APPROXIMATION \\ 
WITH APPLICATIONS TO INFORMATIONAL DISTANCES. II \\
}

\author{S. G. Bobkov$^{1}$}
\thanks{1) 
School of Mathematics, University of Minnesota, USA;
research was partially supported by the Simons Foundation
and the NSF grant DMS-1855575}
\address
{Sergey G. Bobkov \newline
School of Mathematics, University of Minnesota  \newline 
127 Vincent Hall, 206 Church St. S.E., Minneapolis, MN 55455 USA}
\smallskip
\email{bobkov@math.umn.edu}

\author{G. P. Chistyakov$^{2}$}
\thanks{2) 
Faculty of Mathematics, University of Bielefeld, Germany;
research was partially supported by SFB 1283}
\address
{Gennadiy P. Chistyakov\newline
Fakult\"at f\"ur Mathematik, Universit\"at Bielefeld\newline
Postfach 100131, 33501 Bielefeld, Germany}
\smallskip
\email{chistyak@math.uni-bielefeld.de}

\author{F. G\"otze$^{2}$}
\address
{Friedrich G\"otze\newline
Fakult\"at f\"ur Mathematik, Universit\"at Bielefeld\newline
Postfach 100131, 33501 Bielefeld, Germany}
\email{goetze@math.uni-bielefeld.de}

\subjclass
{Primary 60E, 60F} 
\keywords{$\chi^2$-divergence, Relative entropy, Poisson approximation} 

\begin{abstract}
We explore asymptotically optimal bounds for deviations of distributions
of independent Bernoulli random variables from the Poisson limit in terms 
of the Shannon relative entropy and R\'enyi/relative Tsallis distances
(including Pearson's $\chi^2$). This part generalizes the results obtained 
in Part I and removes any constraints on the parameters of the Bernoulli distributions.
\end{abstract}

\maketitle
\markboth{S. G. Bobkov, G. P. Chistyakov and F. G\"otze}{
Relative entropy and $\chi^2$ divergence from the Poisson law}

\def\theequation{\thesection.\arabic{equation}}
\def\E{{\mathbb E}}
\def\R{{\mathbb R}}
\def\C{{\mathbb C}}
\def\P{{\mathbb P}}
\def\Z{{\mathbb Z}}
\def\H{{\rm H}}
\def\Im{{\rm Im}}
\def\Tr{{\rm Tr}}

\def\k{{\kappa}}
\def\M{{\cal M}}
\def\Var{{\rm Var}}
\def\Ent{{\rm Ent}}
\def\O{{\rm Osc}_\mu}

\def\ep{\varepsilon}
\def\phi{\varphi}
\def\F{{\cal F}}
\def\L{{\cal L}}

\def\be{\begin{equation}}
\def\en{\end{equation}}
\def\bee{\begin{eqnarray*}}
\def\ene{\end{eqnarray*}}

\section{{\bf Introduction}}
\setcounter{equation}{0}

\vskip2mm
\noindent
Let $W = X_1 + \dots + X_n$ be the sum of independent random variables 
$X_j$ taking values $1$ and $0$ with respective probabilities 
$p_j$ and $q_j = 1-p_j$. Thus,
\be
w_k = \P\{W = k\} = 
\sum p_1^{\ep_1} q_1^{1-\ep_1} \dots p_n^{\ep_n} q_n^{1-\ep_n}, \qquad
k = 0,1,\dots,n,
\en
where the summation runs over all 0-1 sequences $\ep_1,\dots,\ep_n$ 
such that $\ep_1 + \dots + \ep_n = k$. 

Denote by $Z$ 
a Poisson random variable with parameter $\lambda = p_1 + \dots + p_n$, i.e.,
taking non-negative integer values wih probabilities
\be
v_k = \P\{Z = k\} = \frac{\lambda^k}{k!}\,e^{-\lambda}, \qquad k = 0,1,\dots
\en
It is well known that, if all $p_j$ are small, the distribution of $Z$
approximates the distribution of $W$ in terms of 
the total variation distance $d(W,Z) = \sum_{k=0}^\infty\, |w_k - v_k|$.
In particular, involving the functional $\lambda_2 = p_1^2 + \dots + p_n^2$,
Barbour and Hall \cite{B-H} derived a two sided bound
\be
\frac{1}{32} \min(1,1/\lambda)\,\lambda_2 \, \leq \,
\frac{1}{2}\, d(W,Z) \, \leq \,
\frac{1-e^{-\lambda}}{\lambda}\,\lambda_2.
\en

There is considerable interest as well in the  question of Poisson approximation
for  (stronger) informational distances, including the R\'enyi divergences, 
or equivalently --  the Tsallis relative entropies in their full hierarchy.
Being well-defined in the setting of abstract measure spaces 
(cf. e.g. \cite{E-H}, \cite{B-C-G1}), in the discrete model specified above 
these important quantities are respectively given for any parameter 
$\alpha > 0$ by
$$
D_\alpha = D_\alpha(W||Z) = \frac{1}{\alpha-1}\,\log
\sum_{k=0}^\infty \Big(\frac{w_k}{v_k}\Big)^\alpha \,v_k
$$
and
$$
T_\alpha = T_\alpha(W||Z) = \frac{1}{\alpha-1}\,\bigg[
\sum_{k=0}^\infty 
\Big(\frac{w_k}{v_k}\Big)^\alpha \,v_k - 1\bigg].
$$
The functions $\alpha \rightarrow D_\alpha$ and
$\alpha \rightarrow T_\alpha = \frac{1}{\alpha-1}\,
(\exp\{(\alpha - 1)\,D_\alpha\} - 1)$ are non-decreasing, and in the
particular cases $\alpha=1$ and $\alpha=2$, we deal with the
more familiar relative entropy (Kullback-Leibler distance) and 
the Pearson $\chi^2$-distance
$$
D = D_1 = T_1 = \sum_{k=0}^\infty w_k \log \frac{w_k}{v_k}, 
\qquad
T_2 = \chi^2 = \sum_{k=0}^\infty \frac{(w_k - v_k)^2}{v_k}.
$$
We refer to \cite{Z-H} and \cite{B-C-G2} for historical references 
related to the lower and upper bounds as in (1.3), as well as to recent 
developments towards the problem of bounding of $D$ and $\chi^2$. 
Here, let us only mention a few results in this direction.

In a rather general asymptotic regime (which is typical in applications),
Borisov and Vorozheĭkin \cite{B-V} observed that $\chi^2$ is
approximately $\frac{1}{2}\,\big(\frac{\lambda_2}{\lambda})^2$,
and more precisely,
$$
\lim \frac{\chi^2}{(\lambda_2/\lambda)^2} \, = \, \frac{1}{2} 
\quad {\rm as} \ \ \lambda^6\, \lambda_2 \rightarrow 0.
$$
On the other hand, Harremo\"es, Johnson and 
Kontoyiannis \cite{H-J-K} have recently derived a universal lower 
bound on the relative entropy, 
$D \geq \frac{1}{4}\,\big(\frac{\lambda_2}{\lambda})^2$.
Here, the constant $\frac{1}{4}$ is best possible and is 
asymptotically attained in the case of equal probabilities $p_j$
\cite{H-R}. It is therefore natural to wonder whether or not 
there are two-sided bounds such as
\be
\frac{1}{4}\, \Big(\frac{\lambda_2}{\lambda}\Big)^2  \leq 
D \leq \chi^2 \leq c\,
\Big(\frac{\lambda_2}{\lambda}\Big)^2.
\en 
This turns out to be true in the the case where $\lambda_2/\lambda$ 
is bounded away from 1. Based on orthogonal expansions in 
Charlier polynomials over the Poisson measure and using  
the  Parseval identity in this context, Zacharovas and Hwang \cite{Z-H}
obtained a superior upper bound
\be
\chi^2 \leq 
2\,(\sqrt{e} - 1)^2\, \Big(\frac{\lambda_2}{\lambda}\Big)^2
\Big(1 - \frac{\lambda_2}{\lambda}\Big)^{-3}
\en
(among other similar results for different distances).
Consequently, if for example 
$\frac{\lambda_2}{\lambda} \leq \frac{1}{2}$, 
then (1.4) is fulfilled with $c = 6.74$. 

The upper estimate such as (1.4) also appears as a consequence
of non-uniform bounds which have been recently studued 
in \cite{B-C-G2}. It was shown there that 
$\frac{w_k - v_k}{v_k}$ is of order at most $\lambda_2$
on a large part of the support of the Poisson measure, especially
when $\lambda$ is large. One of the aims of this paper is to extend 
(1.4) modulo absolute constants to the whole range of 
$(\lambda,\lambda_2)$. To formulate results in a compact form, 
let us use the notation $Q_1 \sim Q_2$, whenever two positive quantities 
are related by $c_1 Q_1 \leq Q_2 \leq c_2 Q_1$ with some absolute 
constants $c_j>0$. Introduce the quantity
$$
F = \frac{\max(1,\lambda)}{\max(1,\lambda-\lambda_2)}.
$$
Clearly, $F \geq 1$. 

\vskip5mm
{\bf Theorem 1.1.} {\sl We have
\be
D \sim
\Big(\frac{\lambda_2}{\lambda}\Big)^2\,(1 + \log F), \qquad
\chi^2 \sim
\Big(\frac{\lambda_2}{\lambda}\Big)^2\,\sqrt{F}.
\en
}

%\vskip2mm
If $\frac{\lambda_2}{\lambda}$ is bounded away from 1,  then $F$ 
is bounded, and (1.6) recovers (1.4). A similar conclusion is also true, 
when $\lambda$ is not large, say $\lambda \leq 10$, which is typical 
for applications (note that for such $\lambda$'s, $\frac{\lambda_2}{\lambda}$ 
may be close to 1, and then (1.5) fails to be optimal).
On the other hand, if these two assumptions on $\lambda$ and $\lambda_2$ are 
violated (which we hence forth call the ``degenerate case"), both distances 
are bounded away from zero and can be large, since then
$$
D \, \sim \, 
\log\frac{\lambda}{\max\{1,\lambda-\lambda_2\}}, \qquad
\chi^2 \, \sim \, 
\bigg(\frac{\lambda}{\max\{1,\lambda-\lambda_2\}}\bigg)^{1/2}.
$$
This shows that the lower bound for $D$ in (1.4) may not be
reversed in general.
Indeed, in the extreme case with all $p_j = 1$, we have 
$\lambda_2 = \lambda = n$. Here $\P\{W = n\} = 1$, hence as 
$n \rightarrow \infty$
$$
D = \log \frac{1}{\P\{Z = n\}} = 
\log\Big(\frac{n!}{n^n}\,e^n\Big) \sim \log n,
$$
$$
\chi^2 = \frac{1}{\P\{Z = n\}} - 1 = 
\frac{n!}{n^n}\,e^n - 1 \sim \sqrt{2\pi n}.
$$

As a next step, we employ the non-uniform bounds of \cite{B-C-G2}
to extend (1.4) and (1.6) to all Tsallis entropies.

\vskip5mm
{\bf Theorem 1.2.} {\sl Given $\alpha\,>\, 1$,
\be
T_\alpha \sim
\Big(\frac{\lambda_2}{\lambda}\Big)^2\,F^{\frac{\alpha - 1}{2}}
\en
with involved constants depending on $\alpha$.
In particular, $T_\alpha \leq c_\alpha \chi^2$ as long as
$\frac{\lambda_2}{\lambda} \leq \frac{1}{2}$.
}

\vskip5mm
Let us finally mention one application of Theorem 1.1 
to the problem of the estimation of the difference of entropies
\be
H(W||Z) = H(Z) - H(W),
\en
where $H$ stands for the Shannon entropy, that is,
$$
H(Z) = -\sum_k v_k \log v_k, \qquad H(W) = -\sum_k w_k \log w_k.
$$
The property that $H(W||Z)$ is positive is a consequence
of the assertion, recently proved by Hillion and Johnson \cite{H-J},
that $H(p) \equiv H(W)$ is 
a concave function of the vector $p = (p_1,\dots,p_n)$. Indeed, since
$H(p)$ is invariant under permutations of $p_j$, this entropy
attains its maximum on the simplex $p_j \geq 0$, $p_1 + \dots + p_n = \lambda$
at the point where all coordinates coincide, that is, for
$p_j = \lambda/n$. But in that case, the distribution of $W$ represents 
the binomial law with parameters $n$ and $\lambda/n$ whose entropy 
is dominated by $H(Z)$, as was shown by Harremo\"es \cite{H2001}.

Thus, the difference of entropies 
in this particular model may be viewed as kind of informational 
distance. Sason proposed to bound $H(W||Z)$ for equal $p_j$'s by means of 
the so-called maximal coupling, cf. \cite{S2}. Here, we show that this distance
may be controlled in terms of $\chi^2(W,Z)$, which together
with the upper bound on the Pearson distance as in (1.4)-(1.5) 
leads to the following estimate.

\vskip5mm
{\bf Corollary 1.3.} {\sl With some constants $C_\lambda$ depending 
only on $\lambda$, we have
\be
H(W||Z) \, \leq \, C_\lambda \frac{\lambda_2}{\lambda}.
\en
If $\lambda_2 \leq \frac{1}{2}\,\lambda$, one may take 
$C_\lambda = C \log(2+\lambda)$ with an absolute constant $C$.
}

\vskip5mm
Below, we start with some general
bounds involving the relative entropy and the Pearson distance (Section 2).
In Section 3, we describe several results obtained in \cite{B-C-G2} 
in the non-degenerated case, and employ there
some bounds for the probability function of the Poisson law.
The remaining parts are devoted to the proof of Theorems 1.1 and 1.2
in the degenerate case (Sections 4-10) and of Corollary 1.3
(Section 11). Thus, the paper is structured as follows:

\vskip2mm
1. \ Introduction

2. \ General bounds on relative entropy and $\chi^2$

3. \ Poisson approximation in the non-degenerate case

4. \ Upper bounds on $D$ and $\chi^2$

5. \ Lower bound on $\chi^2$

6. \ Lower bound on $D$

7. \ Proof of Theorem 1.1

8. \ Tsallis versus Vajda-Pearson

9. \ Estimates of Vajda-Pearson distances

10. Proof of Theorem 1.2

11. Difference of entropies

%--------------------------------- Section 2 -------------------------------
\vskip5mm
\section{{\bf General Bounds on Relative Entropy and $\chi^2$}}
\setcounter{equation}{0}

\vskip2mm
\noindent
Before turning to the problem of lower and upper bounds
for the relative entropy and $\chi^2$-distance, we first collect 
several useful general inequalities. If two discrete random elements 
$W$ and $Z$ in a measurable space $\Omega$ take at most countably 
many values $\omega_k \in \Omega$ with probabilities
$w_k = \P\{W = \omega_k\}$ and $v_k = \P\{Z = \omega_k\}$,
the above distances are defined canonically by
$$
D(W||Z) = \sum_{k} w_k \log \frac{w_k}{v_k}, \qquad
\chi^2(W,Z) = \sum_{k} \frac{(w_k - v_k)^2}{v_k}.
$$

\vskip2mm
{\bf Proposition 2.1}. {\sl We have
\be
-\sum_{w_k < v_k} w_k \log\frac{w_k}{v_k} \le 1.
\en
Moreover,
\be
D(W||Z)  \, \ge \, \frac{1}{2}\, \sum_{k} \frac{(w_k-v_k)^2}{\max\{w_k,v_k\}}.
\en
}

\vskip2mm
{\bf Proof}. Using the Taylor formula for the logarithmic function, write
\bee
\sum_{w_k<v_k} w_k \log\frac{w_k}{v_k}
 & = &
\sum_{w_k<v_k}(v_k-(v_k-w_k))\,\log\Big(1-\frac{v_k-w_k}{v_k}\Big) \nonumber \\
 & & \hskip-5mm = \ 
\sum_{w_k<v_k} (w_k-v_k) + \sum_{w_k<v_k} \sum_{m=2}^{\infty}
\frac{1}{m(m-1)}\,\frac{(v_k-w_k)^m}{v_k^{m-1}}.
\ene
Here
$$
\sum_{w_k<v_k} (w_k-v_k) \, = \, 
-\frac{1}{2}\,\sum_{k=0}^\infty |w_k-v_k| \geq -1,
$$
thus proving the first assertion. Similarly, we have a second identity
\bee
\sum_{w_k>v_k}w_k\log\frac{w_k}{v_k}
 & = &
-\sum_{w_k>v_k}w_k\log\frac{v_k}{w_k} \\
 & = &
-\sum_{w_k>v_k}w_k\log\Big(1-\frac{w_k-v_k}{w_k}\Big)\\
 & = &
\sum_{w_k>v_k}(w_k-v_k) + \sum_{w_k>v_k}\, \sum_{m=2}^{\infty}\,
\frac{1}{m} \frac{(w_k-v_k)^m}{w_k^{m-1}}.
\ene
Adding the two identities, we get
$$
\sum_{k}w_k\log\frac{w_k}{v_k} \, \ge \,
\frac{1}{2}\, \sum_{w_k>v_k} \frac{(w_k-v_k)^2}{w_k} +
\frac{1}{2}\, \sum_{w_k<v_k}\frac{(w_k-v_k)^2}{v_k},
$$
which is the desired inequality (2.2). 
\qed

\vskip5mm
{\bf Proposition 2.2.} {\sl Let $W_1$ and $W_2$ be independent, non-negative, 
integer-valued random variables with finite means, and let $Z_1$ and $Z_2$ 
be independent Poisson random variables with 
$\E Z_1 = \E W_1$ and $\E Z_2 = \E W_2$. Then
\be
D(W_1 + W_2||Z_1 + Z_2) \, \le \, D(W_1||Z_1) + D(W_2||Z_2).
\en
In addition,
\be
\chi^2(W_1+W_2,Z_1+Z_2) + 1 \, \le \,
(\chi^2(W_1,Z_1)+1)(\chi^2(W_2,Z_2)+1).
\en
}

\vskip2mm
For the proof, we refer to Johnson \cite{J}, pp.\,133--134. Let 
us only mention that (2.4) is obtained in \cite{J} in the more general form
$$
\sum_{k=0}^{\infty}
\frac{\P\{W_1+W_2=k\}^{\alpha}}{\P(Z_1+Z_2=k\}^{\alpha-1}} \, \le \,
\sum_{k=0}^{\infty}\frac{\P\{W_1=k\}^{\alpha}}{\P\{Z_1=k\}^{\alpha-1}}\,
\sum_{k=0}^{\infty}\frac{\P\{W_2=k\}^{\alpha}}{\P\{Z_2=k\}^{\alpha-1}}
$$
with arbitrary $\alpha \geq 1$, which represents a Poisson analog of weighted
convolution inequalities due to Andersen \cite{A}. Here, for $\alpha = 1$
there is an equality, and comparing the derivatives of both sides at 
this point, we arrive at the relation (2.3).

%--------------------------- Section 3 --------------------------------
\vskip7mm
\section{{\bf \large Poisson Approximation in the Non-Degenerate Case}}
\setcounter{equation}{0}

\vskip2mm
\noindent
Now, we restrict ourselves to the random variables $W = X_1 + \dots + X_n$ 
and $Z$ with distributions described in (1.1)-(1.2). In particular,
\begin{eqnarray}
\P\{Z = 0\} 
 & = &
e^{-(p_1 + \dots + p_n)} \ = \ e^{-\lambda}, \nonumber \\
\P\{W = 0\} 
 & = &
(1-p_1) \dots (1-p_n) \ \leq \ \P\{Z = 0\}.
\end{eqnarray}

The bounds (1.4) follow from the following two assertions proved in \cite{B-C-G2}. 
To compare the lower and upper bounds, we recall the lower bound (1.4) of
Harremo\"es, Johnson and Kontoyiannis \cite{H-J-K}.

\vskip5mm
{\bf Proposition 3.1.} {\sl If $\max_j p_j \leq \frac{1}{2}$, then
$$
\frac{1}{4}\,\Big(\frac{\lambda_2}{\lambda}\Big)^2 \, \leq \,
D(W||Z) \, \leq \, \chi^2(W,Z) \, \leq \, 
C_\lambda\,\Big(\frac{\lambda_2}{\lambda}\Big)^2,
$$
where $C_\lambda$ depends on $\lambda \geq 0$ and is an increasing
continuous function with $C_0 = 2$. In particular,
if $\lambda \leq 1/2$, then
$$
\chi^2(W,Z) \, \leq \, 15\,\Big(\frac{\lambda_2}{\lambda}\Big)^2.
$$
}

\vskip2mm
{\bf Proposition 3.2.} {\sl If $\lambda \geq 1/2$ and
$\lambda_2 \le \kappa \lambda$ with $\kappa \in (0,1)$, then
\be
\frac{1}{4}\,\Big(\frac{\lambda_2}{\lambda}\Big)^2 \, \leq \, 
 D(W||Z) \, \leq \, \chi^2(W,Z) \, \leq \,
c_\kappa\,\Big(\frac{\lambda_2}{\lambda}\Big)^2.
\en
where one may take $c_\kappa = c\,(1-\kappa)^{-3}$ with some absolute
constant, e.g. $c = 7 \cdot 10^6$.
}

\vskip5mm
A natural approach to the Poisson approximation is based on the 
comparison of characteristic functions. Since the random variables
$W$ and $Z$ assume non-negative integer values only, one may equivalently
consider the associated generating functions, similar as in \cite{B-C-G2}.
The generating function for the Poisson law with parameter
$\lambda>0$ is given by
$$
\varphi(w) = \E\,w^Z = \sum_{k=0}^\infty \P\{Z=k\}\, w^k =
e^{\lambda (w-1)} = 
\prod_{j=1}^n\, e^{p_j (w-1)}, 
$$
which is an entire function of the complex variable $w$.
Correspondingly, the generating function for the distribution of
the random variable $W$ is
$$
g(w) = \E\,w^W = \sum_{k=0}^\infty \P\{W=k\}\, w^k = 
\prod_{j=1}^n\, (q_j + p_j w),
$$
which is a polynomial of degree $n$. Hence, the difference between the 
involved probabilities may be expressed with the help of the contour 
integrals by the Cauchy formula
$$
\P\{W=k\} - \P\{Z=k\} = \int_{|w| = r} 
w^{-k}\, (g(w)-\varphi(w))\,d\mu_r(w),
$$
where $\mu_r$ is the uniform probability measure on the circle $|w| = r$ 
of an arbitrary radius $r>0$. This identity for the difference of
probabilities was used in \cite{B-C-G2} in the derivation of 
the upper bound in (3.2), while here the representation
\be
\P\{W=k\} = \int_{|w| = r} w^{-k}\, g(w)\,d\mu_r(w)
\en
will be particularly helpful in the study of the degenerate case.

When estimating the Poisson probabilities 
$$
f(k) = \P\{Z = k\} = \frac{\lambda^k}{k!}\,e^{-\lambda}
$$
for a fixed parameter $\lambda>0$, it is convenient to use the well-known 
Stirling-type two-sided bound:
\be
\sqrt{2\pi}\,k^{k + \frac{1}{2}}\,e^{-k} \leq k! \leq \, 
e\,k^{k + \frac{1}{2}}\,e^{-k} \qquad (k \geq 1).
\en
In particular, it implies the following Gaussian type estimates 
(cf. \cite{B-C-G2})

\vskip5mm
{\bf Lemma 3.3.} {\sl For all $k \geq 1$,
\be
f(k) \, \leq \, \frac{1}{\sqrt{2\pi k}}.
\en
Moreover, if $1 \leq k \leq 2\lambda$, then
\be
\frac{1}{e\sqrt{k}}\,e^{-\frac{(k-\lambda)^2}{\lambda}} \, \leq \, 
f(k) \, \leq \, \frac{1}{\sqrt{2\pi k}}\,e^{-\frac{(k-\lambda)^2}{3\lambda}}.
\en
Here, the lower bound may be improved in the region $k \geq \lambda$ as
\be
f(k) \, \geq \, \frac{1}{e\sqrt{k}}\,e^{-\frac{(k-\lambda)^2}{2\lambda}}.
\en
}

%------------------------------ Section 4 ----------------------
\vskip7mm
\section{{\bf Upper Bounds on $D$ and $\chi^2$}}
\setcounter{equation}{0}

\vskip2mm
\noindent
We now turn to Theorem 1.2 in the degenerate case, where the optimal
bounds on the relative entropy and $\chi^2$ have a different behavior.
As an intermediate step, let us derive the following upper bounds for the
$\chi^2$-distance and the relative entropy, by using the quantity
$$
Q \, = \, \lambda/\max\{1,\lambda-\lambda_2\}.
$$

\vskip5mm
{\bf Proposition 4.1}. {\sl For $\lambda\ge 1/2$, we have
\begin{eqnarray}
\chi^2(W,Z) 
 & \le & 
19\sqrt{Q}, \\
D(W||Z) 
 & \le & 
23\,\log(e Q).
\end{eqnarray}
}

\vskip2mm
These bounds are sharp when $\lambda_2 \geq \kappa \lambda$, 
cf. Propositions 5.1 and 6.1.

\vskip5mm
{\bf Proof}.
Setting $g(w) = \prod_{l=1}^n (q_l + p_l w)$, $w \in \C$, we exploit the
contour integral representation (3.3), i.e.,
$$
\P\{W=k\} = 
\frac{1}{2\pi}\, r^{-k}\int_{-\pi}^\pi g(re^{i\theta})\,e^{-ik\theta}\,d\theta,
\qquad r>0.
$$
It yields an upper bound
\be
\P\{W=k\} \leq R_k(r)\,I(r),
\en
where
$$
R_k(r) = r^{-k} \prod_{l=1}^n(q_l+p_lr) \quad {\rm and} \quad
I(r) = \frac{1}{2\pi}\int_{-\pi}^{\pi}
\prod_{l=1}^n \frac{|q_l + p_l r e^{i\theta}|}{q_l + p_l r}\,d\theta.
$$
Let us choose $r=k/\lambda$. Since $q_j+p_j r\le e^{p_j(r-1)}$, we have
$$
R_k(r)\le r^{-k}\prod_{j=1}^n(q_j+p_j r)\le e^{\lambda(r-1)-k\log r}=
\Big(\frac{e\lambda}{k}\Big)^k\,e^{-\lambda}.
$$
Moreover, applying $(\frac{e}{k})^k \leq e\sqrt{k}\,\frac{1}{k!}$, 
cf. (3.4), the above is simplified to
\be
R_k(r)\le e\sqrt{k}\ \frac{\lambda^k}{k!}\,e^{-\lambda} =
e\sqrt{k}\, f(k),
\en
where $f(k)$ is the density of the Poisson law with parameter $\lambda$.

Now, to bound $I(r)$, for all $|\theta| \leq \pi$, using
$\sin(\frac{\theta}{2}) \geq \frac{1}{\pi}\,\theta$, we have 
\bee
\prod_{l=1}^n \frac{|q_l + p_l r\, e^{i\theta}|}{q_l+p_lr}
 & = &
\prod_{l=1}^n 
\Big(1-\frac{4q_l p_l r}{(q_l+p_lr)^2}\, \sin^2\frac{\theta}{2}\Big)^{1/2} \\
 & \le &
\exp\Big\{-2\,\sin^2\frac{\theta}{2}\,\sum_{l=1}^n
\frac{q_l p_l\, r}{(q_l + p_l r)^2}\Big\} 
 \, \le \,
\exp\Big\{-\frac{2\theta^2}{\pi^2} \sum_{l=1}^n
\frac{q_l p_l\, r}{(q_l + p_l r)^2}\Big\}.
\ene
Here
$$
\sum_{l=1}^n \frac{q_l p_l\, r}{(q_l + p_l r)^2} \, \geq \,
\frac{1}{r}\, \sum_{l=1}^n q_l p_l = \frac{1}{r}\, (\lambda - \lambda_2) \quad
{\rm in \ case} \ \ r \geq 1
$$
and
$$
\sum_{l=1}^n \frac{q_l p_l\, r}{(q_l + p_l r)^2} \, \geq \,
r\, \sum_{l=1}^n q_l p_l = r\, (\lambda - \lambda_2) \quad
{\rm in \ case} \ \ r \leq 1.
$$
These right-hand sides have the form
$$
\psi(r) = \min\{r,1/r\}\,(\lambda - \lambda_2),
$$ 
and we get
\bee
I(r) 
 & \le & 
\frac{1}{2\pi} \int_{-\pi}^{\pi}
\exp\Big\{-\frac{2}{\pi^2}\,\psi(r)\,\theta^2\Big\}\,d\theta
 \, = \,
\frac{1}{4\,\psi(r)^{1/2}} \int_{-2\sqrt{\psi(r)}}^{2\sqrt{\psi(r)}}
e^{-\frac{1}{2}\,x^2}\,dx \\
 & \leq &
\frac{1}{4\,\psi(r)^{1/2}}\, \min\big\{\sqrt{2\pi},4\,\psi(r)^{1/2}\big\}
 \, \leq \,
\min\big\{1,\psi(r)^{-1/2}\big\}.
\ene

First, we consider the region $\frac{1}{4}\,\lambda \le k\le 4\lambda$,
in which case $\frac{1}{4} \le r \le 4$ and 
$\psi(r) \geq \frac{1}{4}\,(\lambda - \lambda_2)$ and thus
$$
I(r) \leq \min\Big\{1,\frac{2}{\sqrt{\lambda - \lambda_2}}\Big\}
 \, \le \,2\sqrt{Q_0}, \qquad
Q_0 = 1/\max\{1,\lambda-\lambda_2\}.
$$ 
Applying this bound together with (4.4) in (4.3), we get
\be
\P\{W=k\} \, \le \, 2e\sqrt{Q_0} \ \sqrt{k}\, f(k).
\en

As for the regions $1 \le k < \frac{1}{4}\,\lambda$ and
$k>4\lambda$, we use the property $|I(r)|\le 1$, which yields 
simpler upper bounds
\be
\P\{W=k\} \, \le \,\Big(\frac{e\lambda}{k}\Big)^k\,e^{-\lambda}
 \, \leq \, e\sqrt{k}\, f(k).
\en

Now, recall that $\P\{W=0\} \leq f(0)$ (as mentioned in (3.1)) and write
\bee
\chi^2(W,Z)
 & = &
\sum_{k=0}^{\infty}\frac{\P\{W=k\}^2}{f(k)}-1 \, \leq \, S_1+S_2+S_3 \\
 & = &
\bigg(\sum_{1 \le k < \frac{1}{4}\,\lambda} +
\sum_{\frac{1}{4}\,\lambda \leq k \leq 4\lambda} +
\sum_{k > 4\lambda}\bigg)\, \frac{\P\{W=k\}^2}{f(k)}.
\ene
By (4.5),
$$
S_2 \, \le \, 2e\sqrt{Q_0}
\sum_{\frac{1}{4}\,\lambda \leq k \leq 4\lambda} \sqrt{k}\ \P\{W=k\}
 \, \le \,
4e\sqrt{Q} \sum_{\frac{1}{4}\,\lambda \leq k \leq 4\lambda} \P\{W=k\} \, \leq \, 
4e\sqrt{Q}.
$$

To estimate $S_1$, first note that $S_1 = 0$ for $\lambda < 4$.
For $\lambda \geq 4$, using the property that the function 
$k \rightarrow (\frac{e\lambda}{k})^k$ is increasing for $k < \lambda$,
we obtain from (4.6) that
\bee
S_1
 & \le & 
e^{-\lambda+1} \sum_{k < \frac{1}{4}\,\lambda}
\sqrt{k}\,\Big(\frac{e\lambda}k\Big)^k \, \le \,
\frac{1}{2}\sqrt{\lambda}\ e^{-\lambda+1}
\sum_{1 \leq k < \frac{1}{4}\,\lambda} \Big(\frac{e\lambda}{k}\Big)^k \\
 & \leq &
\frac{1}{2}\sqrt{\lambda}\ e^{-\lambda+1}
\sum_{1 \leq k < \frac{1}{4}\,\lambda} (4e)^{\lambda/4}  \, \leq \, 
e\,\Big(\frac{\lambda}{4}\Big)^{3/2}\,\Big(\frac{4}{e^3}\Big)^{\lambda/4} \\
 & \leq &
e\,\Big(\frac{3}{2e\log(e^3/4)}\Big)^{3/2} \, < \, 0.544.
\ene
Here we applied the inequality
\be
x^p\, c^{x} \, \leq \, \Big(\frac{p}{e\log(1/c)}\Big)^p \qquad p,x>0, \ \ 
0 < c < 1,
\en
with $p = 3/2$ and $c = 4/e^3$.

To estimate $S_3$, one may bound the sequence $\sqrt{k}\,(\frac{e\lambda}{k})^k$
for $k > 4\lambda \geq 2$ by the geometric progression $A b^k$ with
suitable parameters $A>0$ and $0 < b < 1$.
To this aim, consider the function
\bee
u(x)
 & = & 
\log\Big(\sqrt{x}\,\Big(\frac{e\lambda}{x}\Big)^x\Big) - \log(b^x) \\
 & = &
\frac{1}{2}\,\log x + x + x\log \lambda - x\log x - x\log b, \qquad x \geq 4\lambda.
\ene
We have
$$
u'(x) \, = \, \frac{1}{2x} + \log \lambda - \log x - \log b \, \leq \,
\frac{1}{4} + \log\frac{1}{4b} \, \leq \, 0,
$$
if $b \geq \frac{1}{4}\,e^{1/4}$ which we assume. In this case, 
$u$ is decreasing, so that
$u(x) \leq u(4\lambda) = 
\log\big(2\sqrt{\lambda}\,(\frac{e}{4b})^{4\lambda}\big) \leq \log A$,
where
$$
A \, = \, 2 \sup_{\lambda \geq 1/2} 
\sqrt{\lambda}\,\Big(\frac{e}{4b}\Big)^{4\lambda} \, = \,
\sup_{y \geq 2} \sqrt{y}\,\Big(\frac{e}{4b}\Big)^y \, = \,
\Big(\frac{1}{2e\log(3/e)}\Big)^{1/2} \, < \, 1.366,
$$
where on the last step we choose $b = 3/4$ and applied (4.7)
with $p = 1/2$ and $c = e/3$. Thus, putting $k_0 = [4\lambda] + 1$
and noting that $k_0 \geq 2$, we get
\bee
S_3 
 & \le &
e^{-\lambda+1}\sum_{k > 4\lambda} \sqrt{k}\,
\Big(\frac{e\lambda}{k}\Big)^k \, \le \,
\sqrt{e} \sum_{k \geq k_0} A\,\Big(\frac{3}{4}\Big)^k \\
 & = &
4 A\sqrt{e} \,\Big(\frac{3}{4}\Big)^{k_0} \, \leq \, 
\frac{9}{4}\,A \sqrt{e} \, < \, 5.067.
\ene

Finally, using $Q = \lambda Q_0 \geq 1/2$ 
(due to $\lambda \geq 1/2$), we get $S_1 + S_3 < 5.611 \leq 5.611\sqrt{2Q}$.
This gives $S_1 + S_2 + S_3 < (5.611\sqrt{2} + 4e)\sqrt{Q} < 18.81\sqrt{Q}$, 
so (4.1) follows.

Turning to the second assertion and using $\P\{W=0\} \leq f(0)$, 
write similarly
\bee
D(W||Z)
 & = &
\sum_{k=0}^{\infty} \P\{W=k\}\, \log\frac{\P\{W=k\}}{\P\{Z=k\}}
 \, = \, T_1 + T_2 + T_3\\
 & \leq &
\Big(\sum_{1 \leq k < \frac{1}{4}\,\lambda} +
\sum_{\frac{1}{4}\,\lambda \leq k \leq 4\lambda} +
\sum_{k > 4\lambda}\Big)\ \P\{W=k\}\,\log\frac{\P\{W=k\}}{f(k)}.
\ene
For the region $\frac{1}{4}\,\lambda \le k\le 4\lambda$, we can
apply the bound (4.5) again, which gives
$$
\P\{W=k\} \, \le \, 2\sqrt{Q_0} \
e\sqrt{k}\, f(k)  \, \le \, 4e\sqrt{Q}\, f(k),
$$
and therefore, using $Q \geq 1/2$, 
$$
T_2 \, \leq \, \log(4e) + \frac{1}{2}\,\log Q \, \leq \, 
\frac{\log(4e) - \frac{1}{2}\,\log 2}{\log(e/2)}\, \log(eQ) \, < \,
6.65\,\log(eQ).
$$
Using (4.6) together with the inequality $\log(et) \leq t$
($t > 0$), we obtain, similarly to the derivation of the bound
on $T_1$ in the $\chi^2$-case, that 
\bee
T_1
 & \le &
e^{-\lambda}\, \sum_{1 \leq k < \frac{1}{4}\,\lambda}
\Big(\frac{e\lambda}{k}\Big)^k\,\log(e\sqrt{k})
 \, \le \,
e^{-\lambda} \log(e\sqrt{\lambda/4})\
\sum_{1 \leq k < \frac{1}{4}\,\lambda} \Big(\frac{e\lambda}{k}\Big)^k \\
 & \le &
\Big(\frac{\lambda}{4}\Big)^{3/2}\,
\,\Big(\frac{4}{e^3}\Big)^{\lambda/4} \, \leq \, 
\Big(\frac{3}{2e\log(e^3/4)}\Big)^{3/2} \, < \, 0.2.
\ene
Choosing again $k_0 = [4\lambda] + 1$ similarly to
the derivation of the bound on $S_3$ in the $\chi^2$-case, we also get
\bee
T_3 
 & \le &
e^{-\lambda}\, \sum_{k > 4\lambda}
\Big(\frac{e\lambda}{k}\Big)^k\,\log(e\sqrt{k}) \leq 
e^{-\lambda+1}\sum_{k \geq k_0} \sqrt{k}\,
\Big(\frac{e\lambda}{k}\Big)^k \, < \, 5.067.
\ene
Hence,
$T_1 + T_3 < 5.087 < 16.578\,\log(eQ)$, and (4.2) follows as well.
\qed

%----------------------------- Section 5 ------------------------------------
\vskip7mm
\section{{\bf \large Lower Bound on $\chi^2$}}
\setcounter{equation}{0}

\vskip2mm
\noindent
Here, we complement Proposition 4.1 by a similar lower bound for the
$\chi^2$-distance in terms of the same quantity
$Q = \lambda/\max\{1,\lambda-\lambda_2\}$.
Let $c_0 = 2.5 \cdot 10^{-6}$.

\vskip5mm
{\bf Proposition 5.1}. {\sl
If $\lambda\ge 1/2$, then with some absolute constant 
$c \in [c_0,1)$
\be
1 + \chi^2(W,Z)\ge c\sqrt{Q}.
\en
Moreover,
\be
\chi^2(W,Z) \ge \frac{c}{9} \sqrt{Q}
\en
as long as $\lambda_2 \geq (1 - \frac{c^2}{4})\, \lambda$.
}

\vskip5mm
Suppose that $\lambda_2 \geq (1 - \frac{c^2}{4})\, \lambda$. To derive 
(5.2) from (5.1), it is sufficient to require that $c\sqrt{Q} \geq 2$, 
since then $c\sqrt{Q} - 1 \geq \frac{c}{2}\sqrt{Q}$. This condition is fulfilled, 
as long as $\lambda \geq \lambda_0 = \frac{4}{c^2}$ and then we obtain (5.2).
In the remaining case $\frac{1}{2} \leq \lambda \leq \lambda_0$, 
the inequality (5.2) follows from the lower bound
$$
\chi^2(W,Z) \ge \frac{1}{4}\,\Big(\frac{\lambda_2}{\lambda}\Big)^2,
$$ 
cf. (1.4). Indeed, in this case, 
$\lambda - \lambda_2 \leq \frac{c^2}{4}\,\lambda \leq 1$, so that $
Q = \lambda \leq \frac{4}{c^2}$, and thus
$\frac{c}{9} \sqrt{Q} \leq \frac{2}{9}$, while
$\frac{1}{4}\,(\frac{\lambda_2}{\lambda})^2 \geq \frac{1}{4}\,(1 - \frac{c^2}{4})^2$.

\vskip2mm
Thus, it remains to derive the first inequality (5.1). First we shall prove it, 
assuming that $\lambda-\lambda_2$ is sufficiently large. As in Section 4, 
for any fixed $r>0$, we apply the Cauchy theorem and write
$$
\P\{W=k\} = \int_{|w|=r} w^{-k} \prod_{l=1}^n(q_l + p_l w)\,d\mu_r(w) = 
R_k(r)\,I_k(r)
$$
with integration over the uniform distribution $\mu_r$ on the circle $|w|=r$
of the complex plane. Here and below
$$
R_k(r)=r^{-k}\prod_{l=1}^n(q_l+p_l r)
$$
and
$$
I_k(r) \, = \,
\frac{1}{2\pi} \int_{-\pi}^{\pi} \prod_{l=1}^n
\frac{|q_l + p_l r\, e^{i\theta}|}{q_l + p_l r}\,
\exp\Big\{-ik\theta + i\sum_{l=1}^n 
{\rm Im}\big(\log(q_l + p_l r\, e^{i\theta})\big)\Big\}\,d\theta.
$$
We split the integration over the two regions so that to work with the
representation
$$
\P\{W=k\} = R_k(r)\,I_k(r) = R_k(r)\,(I_{k1}(r)+I_{k2}(r)),
$$
where 
\bee
I_{k1}(r)
 & \hskip-2mm = &  \hskip-2mm
\frac{1}{2\pi} \int_{-\frac{\pi}{2}}^{\frac{\pi}{2}}\, \prod_{l=1}^n
\frac{|q_l + p_l r\, e^{i\theta}|}{q_l + p_l r}\,
\exp\Big\{-ik\theta + i\sum_{l=1}^n 
{\rm Im}\big(\log(q_l + p_l r\, e^{i\theta})\big)\Big\}\,d\theta,\\
I_{k2}(r)
 & \hskip-2mm = &  \hskip-2mm
\frac{1}{2\pi} \int_{\frac{\pi}{2} < |\theta| < \pi}\, \prod_{l=1}^n
\frac{|q_l + p_l r\, e^{i\theta}|}{q_l + p_l r}\,
\exp\Big\{-ik\theta + i\sum_{l=1}^n 
{\rm Im}\big(\log(q_l + p_l r e^{i\theta})\big)\Big\}\, d\theta.
\ene
To properly estimate $I_k(r)$ from below,
$I_{k2}(r)$ needs to be estimated from above (in absolute value), while
$I_{k1}(r)$, which is a real number, should be estimated from below.

Furthermore, the quantity $R_k(r)$ needs to be estimated from below
as well. To this aim, we choose 
the radius $r = r(k) > 0$ by the condition $R_k'(r)=0$, or equivalently
\be
F(r) \equiv \sum_{l=1}^n\frac{p_l r}{q_l + p_l r} = k.
\en
Since the function $F$ is monotone and $F(0)=0$, $F(\infty)= n$,
there is a unique solution, say $r$, to this equation as long as $n>k$ 
(which may be assumed). We also assume that not all $p_k$ are 
equal to 0 or 1, so that $\lambda_2 < \lambda$.

Let us also emphasize that $F$ is concave on the positive half-axis. 
Since $F(1) = \lambda$, we necessarily have $r(k) < 1$ in case $k < \lambda$, 
and $r(k) > 1$ in case $k > \lambda$.

\vskip5mm
{\bf Lemma 5.2.} {\sl For any $k = 0,\dots,n-1$, 
the solution $r=r(k)$ to the equation $(5.3)$ satisfies 
$$
r \geq 1 + \frac{k-\lambda}{\lambda-\lambda_2}.
$$
Moreover, in case $|k-\lambda| \le \frac{1}{6}\, (\lambda-\lambda_2)$, 
we have $\frac{5}{6} \leq r \leq \frac{6}{5}$, and actually with some 
$0 \leq b_i \leq 1$
\bee
r 
 & = &
1 + \Big(\frac{6}{5}\Big)^2\,b_1\,\frac{k-\lambda}{\lambda-\lambda_2} \\
 & = &
1 + \frac{k-\lambda}{\lambda-\lambda_2} + 
\Big(\frac{6}{5}\Big)^9\,b_2\, \frac{\lambda_2-\lambda_3}{\lambda-\lambda_2}\,
\Big(\frac{k-\lambda}{\lambda-\lambda_2}\Big)^2.
\ene
}

%\vskip5mm
{\bf Proof.} We have
$$
F'(r) = \sum_{l=1}^n\frac{p_l q_l}{(q_l + p_l r)^2}, \qquad 
F'(1) = \lambda - \lambda_2.
$$
The inverse function $F^{-1}:[0,n) \rightarrow [0,\infty)$ is
increasing and convex. Hence, for any $s \in [0,n)$,
\bee
F^{-1}(s) 
 & \geq &
F^{-1}(\lambda) + (F^{-1})'(\lambda)\,(s - \lambda) \\
 & = &
F^{-1}(\lambda) + \frac{1}{F'(F^{-1}(\lambda))}\,(s - \lambda) \, = \,
1 + \frac{1}{\lambda - \lambda_2}\,(s - \lambda).
\ene
Plugging $s = k$, we obtain the first inequality.

Now, since $q_l + p_l r \leq 1$ for $r \leq 1$, we conclude that 
$F'(r) \geq \sum_{l=1}^n p_l q_l = \lambda - \lambda_2$ 
and $F(1) - F(r) \geq (1-r)(\lambda - \lambda_2)$. Thus, if $k \leq \lambda$, 
we obtain that
$$
\frac{1}{6}\, (\lambda-\lambda_2) \geq |k-\lambda| = F(1) - F(r(k)) 
\geq (1-r(k))(\lambda - \lambda_2),
$$
implying $r(k) \geq \frac{5}{6}$. For $r \geq 1$, one may use $q_l + p_l r \leq r$,
which gives $F'(r) \geq \frac{1}{r^2}\,(\lambda - \lambda_2)$ and
$F(r) - F(1) \geq (1 - \frac{1}{r})\,(\lambda - \lambda_2)$.
Hence, again by the assumption,
$$
\frac{1}{6}\, (\lambda-\lambda_2) \geq k-\lambda = F(r(k)) - F(1) 
\geq \Big(1 - \frac{1}{r(k)}\Big)\,(\lambda - \lambda_2),
$$
implying $r(k) \leq \frac{6}{5}$. In both cases, 
$\frac{5}{6} \leq r(k) \leq \frac{6}{5}$, proving the second assertion of the lemma.

Now, in the interval $\frac{5}{6} \leq r \leq \frac{6}{5}$, we necessarily have
$\frac{5}{6} \leq q_l + p_l r \leq \frac{6}{5}$, so that
$$
\Big(\frac{5}{6}\Big)^2\,(\lambda - \lambda_2) \leq F'(r) \leq 
\Big(\frac{6}{5}\Big)^2\,(\lambda - \lambda_2). 
$$
In addition,
$$
-F''(r) \, = \, 
2\sum_{l=1}^n\frac{p_l^2 q_l}{(q_l + p_l r)^3} \, \leq \,
2 \cdot \Big(\frac{6}{5}\Big)^3 \sum_{l=1}^n p_l^2 q_l \, = \, 
2 \cdot \Big(\frac{6}{5}\Big)^3\,(\lambda_2 - \lambda_3).
$$

Let us now write the Taylor expansion up to the linear and quadratic terms 
for the inverse function $F^{-1}(s)$ around the point $\lambda$. Then we get
\bee
F^{-1}(s) 
 & = & 
1 + \frac{1}{F'(F^{-1}(s_1))}\,(s - \lambda) \\
 & = &
1 + \frac{1}{F'(1)}\,(s - \lambda) - 
\frac{1}{2\,F'(F^{-1}(s_2))^3}\,F''(F^{-1}(s_2))\,(s - \lambda)^2,
\ene
where the points $s_1$ and $s_2$ lie between $\lambda$ and $s$. 
Putting $r = F^{-1}(s)$ and $r_i = F^{-1}(s_i)$, the above is simplified as
\bee
r 
 & = &
1 + \frac{1}{F'(r_1)}\,(s - \lambda) \\
 & = &
1 + \frac{1}{\lambda - \lambda_2}\,(s - \lambda) - 
\frac{1}{2 F'(r_2)^3}\,F''(r_2)\,(s - \lambda)^2
\ene
where $r_1$ and $r_2$ lie between $1$ and $r$. It remains to apply these 
equalities with $s = k$, that is, $r = r(k)$, and note that
$\frac{1}{F'(r_1)} \leq (\frac{6}{5})^2\,\frac{1}{\lambda - \lambda_2}$,
while
$$
\frac{1}{2\,F'(r_2)^3}\,|F''(r_2)| \, \leq \, 
\frac{1}{2\,(\frac{5}{6})^6\,(\lambda - \lambda_2)^3}\cdot
2 \cdot \Big(\frac{6}{5}\Big)^3\,(\lambda_2 - \lambda_3) \, = \, 
\Big(\frac{6}{5}\Big)^9\,
\frac{\lambda_2 - \lambda_3}{(\lambda - \lambda_2)^3}.
$$
Note that $(\frac{6}{5})^2 = 1.44$ and $(\frac{6}{5})^9 < 5.16$.
\qed

\vskip5mm
{\bf Lemma 5.3.} {\sl Let $r=r(k)$ be the solution of $(5.3)$
for $0 \le \lambda - k \le \frac{1}{6}\,(\lambda-\lambda_2)$. Then
$$
R_k(r) \, = \, r^{-k}\prod_{l=1}^n(q_l+p_l r) \, \ge \, 
\exp\Big\{-4\,\frac{(\lambda-k)^2}{\lambda-\lambda_2}\Big\}.
$$
}

%\vskip2mm
{\bf Proof.} The function
$$
\psi_k(r) = \log R_k(r) = \sum_{l=1}^n \log(q_l+p_l r) - k\log r, \qquad
r > 0,
$$
is vanishing at $r=1$ and has derivative
$$
\psi_k'(r) = \sum_{l=1}^n \frac{p_l}{q_l+p_l r} - \frac{k}{r} = 
\frac{F(r) - k}{r} = \frac{F(r) - F(r(k))}{r}.
$$
Since $F$ is increasing and concave, $F(a) - F(b) \leq F'(b)\,(a-b)$
whenever $a \geq b > 0$. In particular, in the interval $r(k) \leq r \leq 1$,
we have 
$$
\psi_k'(r) \leq \frac{F'(r(k))}{r}\,(r - r(k)) \leq 
\frac{F'(r(k))}{r(k)}\,(1 - r(k)),
$$
which implies
$$
\psi_k(r(k)) = \psi_k(r(k)) - \psi_k(1) \geq 
-\frac{F'(r(k))}{r(k)}\,(1 - r(k))^2.
$$
By Lemma 5.2, $\frac{5}{6} \leq r(k) \leq 1$ and 
$1 - r(k) \leq (\frac{6}{5})^2\,\frac{\lambda - k}{\lambda-\lambda_2}$.
Moreover, as was shown in the proof, 
$F'(r(k)) \leq (\frac{6}{5})^2\,(\lambda - \lambda_2)$. Hence
$$
\frac{F'(r(k))}{r(k)}\,(1 - r(k))^2 \leq 
\frac{(\frac{6}{5})^2\,(\lambda - \lambda_2)}{5/6}\,
\Big(\Big(\frac{6}{5}\Big)^2\,\frac{k-\lambda}{\lambda-\lambda_2}\Big)^2 =
\Big(\frac{6}{5}\Big)^7\,\frac{(k-\lambda)^2}{\lambda-\lambda_2}.
$$
Here, $(\frac{6}{5})^7 < 3.6$.
\qed

\vskip5mm
{\bf Lemma 5.4.} {\sl 
Let $\lambda-\lambda_2 \ge 100$. Then, for 
$0 \le \lambda-k \le \frac{1}{6}\,(\lambda-\lambda_2)$,
$$
I_k(r(k)) \ge \frac{1}{10\sqrt{\lambda-\lambda_2}}.
$$
}

%\vskip5mm
{\bf Proof.}
By Lemma 5.2, $1 \geq r(k) \geq \frac{5}{6}$. As in the proof of 
Proposition 4.1, recall that for $r>0$ and $-\pi\le\theta\le\pi$,
\bee
\prod_{l=1}^n \frac{|q_l + p_l r\,e^{i\theta}|}{q_l+p_lr}
 & = &
\prod_{l=1}^n
\Big(1 - \frac{4q_l p_l r}{(q_l + p_l r)^2}\,\sin^2\frac{\theta}2\Big)^{1/2} \\
 & \le & 
\exp\Big\{-2\sum_{l=1}^n
\frac{q_l p_l r}{(q_l + p_l r)^2}\,\sin^2\frac{\theta}{2}\,\Big\}.
\ene
For $\frac{5}{6} \le r\le 1$, necessarily $q_l + p_l r \le 1$ and
$$
\sum_{l=1}^n \frac{q_l p_l r}{(q_l+p_lr)^2} \, \ge \,
\sum_{l=1}^n q_l p_l r \, = \,
(\lambda-\lambda_2)\, r \, \ge \, \frac{5}{6}\,(\lambda-\lambda_2).
$$
Hence
\bee
|I_{k2}(r)|
 & \le &
\frac{1}{2\pi} \int_{\frac{\pi}{2} \le |\theta| \le \pi} \,
\prod_{l=1}^n\frac{|q_l + p_l r\,e^{i\theta}|}{q_l + p_l\,r}\,d\theta \\
 & \le &
\frac{1}{2\pi} \int_{\frac{\pi}{2} \le |\theta|\le \pi}
\exp\Big\{-\frac{5}{3}\,(\lambda-\lambda_2)\,\sin^2\frac{\theta}{2}\Big\}\,d\theta 
 \, \le \,
\frac{1}{2}\, e^{-\frac{5}{6}\,(\lambda-\lambda_2)}.
\ene

Let us now estimate $I_{k1}(r)$ from below. 
Using $4q_l p_l r \leq (q_l+p_lr)^2$
which is the same as $(q_l - p_l r)^2 \geq 0$, we have, for $|\theta| \le \pi/2$,
$$
\frac{4q_lp_l r}{(q_l+p_lr)^2}\, \sin^2\frac{\theta}{2} \, \le \,
\frac{1}{2}, \quad l=1,\dots,n.
$$
In the region $0 \leq \ep \leq \ep_0 < 1$, there is a lower bound 
$1-\ep \geq e^{-c\ep}$ with best attainable constant when $\ep = \ep_0$.
In the case $\ep_0 = \frac{1}{2}$, this constant is given by
$c = 2\log 2$. Therefore, for $|\theta|\le\frac{\pi}{2}$,
$$
\prod_{l=1}^n \frac{|q_l +p_l r\,e^{i\theta}|}{q_l + p_l r}
 \, \ge \,
\exp\Big\{-\log 2\,
\sum_{l=1}^n \frac{4q_lp_l r}{(q_l + p_l r)^2}\,\sin^2(\theta/2)\Big\}.
$$
Here, the involved function 
$$
w_l(r) = \frac{r}{(q_l + p_l r)^2}, \qquad r \geq 0,
$$ 
is increasing in 
$0 \leq r \leq r_l \equiv q_l/p_l$ and decreasing in $r \geq r_l$. Hence,
if $r_l \geq 1$, then $\max_{\frac{5}{6} \leq r \leq 1} w_l(r) = w_l(1) = 1$.
If $r_l \leq \frac{5}{6}$, that is, when $p_l \geq \frac{6}{11}$,
we have
$$
\max_{\frac{5}{6} \leq r \leq 1} w_l(r) = w_l(5/6) = 
\frac{\frac{5}{6}}{(q_l + p_l \, \frac{5}{6})^2} \leq \frac{6}{5}.
$$
Finally, if $\frac{5}{6} \leq r_l \leq 1$, which is equivalent to 
$ \frac{1}{2} \leq p_l \leq \frac{6}{11}$, we have
$$
\max_{\frac{5}{6} \leq r \leq 1} w_l(r) = w_l(r_l) = 
\frac{1}{4p_l q_l} \leq
\frac{1}{4 \cdot \frac{6}{11} \cdot \frac{5}{11}} = \frac{121}{120}.
$$
Thus, in all cases, $w_l(r) \leq \frac{6}{5}$ on the interval 
$\frac{5}{6} \leq r \leq 1$, so that
\bee
\prod_{l=1}^n \frac{|q_l +p_l r\,e^{i\theta}|}{q_l + p_l r} 
 & \ge &
\exp\Big\{-\frac{6}{5} \log 2 \sum_{l=1}^n 4q_l p_l\,\sin^2(\theta/2)\Big\} \\
 & \ge &
\exp\Big\{-\frac{6}{5}\, (\log 2)\,(\lambda-\lambda_2)\,\theta^2\Big\},
\ene
and thus
\bee
\frac{1}{2\pi} \int_{-\frac{\pi}{2}}^{\frac{\pi}{2}} \, \prod_{l=1}^n
\frac{|q_l + p_l r\,e^{i\theta}|}{q_l + p_l r}\,d\theta
 & \ge &
\frac{1}{2\pi} \int_{-\frac{\pi}{2}}^{\frac{\pi}{2}}
\exp\Big\{-\frac{6}{5}\,(\log 2)\,(\lambda-\lambda_2)\,\theta^2\Big\}\,d\theta\\
 & & \hskip-32mm = \ 
\frac{1}{2\pi\sqrt{\frac{6}{5}\,(\log 4)\,(\lambda-\lambda_2)}} \, 
\int_{-\frac{\pi}{2} \sqrt{\frac{6}{5}\,(\log 4)\,
(\lambda-\lambda_2)}}^{\frac{\pi}{2} \sqrt{\frac{6}{5}\,(\log 4)\,
(\lambda-\lambda_2)}} \exp\Big\{-\frac{1}{2}\,x^2\Big\}\,dx \\
 & & \hskip-32mm \geq \ 
0.3093\, \frac{1}{\sqrt{\lambda-\lambda_2}}.
\ene
Here we used $\lambda - \lambda_2 \geq 100$, which ensures that
\bee
\frac{1}{2\pi\sqrt{\frac{6}{5}\,\log 4}}
\int_{-\frac{\pi}{2}\sqrt{\frac{6}{5}\,(\log 4)\,
(\lambda-\lambda_2)}}^{\frac{\pi}{2}\sqrt{\frac{6}{5}\,(\log 4)\,
(\lambda-\lambda_2)}} \ e^{-\frac{1}{2}\,x^2}\,dx
 & \geq & 
\frac{1}{2\pi\sqrt{\frac{6}{5}\,\log 4}}
\int_{-5\pi \sqrt{\frac{6}{5}\,\log 4}}^{5\pi \sqrt{\frac{6}{5}\,\log 4}}
e^{-\frac{1}{2}\,x^2}\,dx \\
 & & \hskip-30mm = \ 
\frac{1}{\sqrt{2\pi \, \frac{6}{5}\, \log 4}}\
\P\Big\{|\xi| \leq 5\pi \sqrt{\frac{6}{5}\,\log 4}\Big\} \, > \, 0.3093,
\ene
where $\xi \sim N(0,1)$. In addition (recalling one of the upper
bounds when bounding the integral $I_{k2}$ from above), and using
$\sin(\theta/2) \geq \frac{\sqrt{2}}{\pi}\,\theta$ for $0 \leq \theta \leq \pi/2$,
we get that
\bee
\frac{1}{\pi} \int_{-\frac{\pi}{2}}^{\frac{\pi}{2}} \, \prod_{l=1}^n
\frac{|q_l + p_l r\,e^{i\theta}|}{q_l + p_l r}\ \theta^6\,d\theta
 & \le &
\frac{1}{\pi} \int_{-\frac{\pi}{2}}^{\frac{\pi}{2}}
\exp\Big\{-\frac{5}{3}\,(\lambda-\lambda_2)\,\sin^2\frac{\theta}{2}\Big\}\ 
\theta^6\,d\theta \\
 & \le &
\frac{1}{\pi} \int_{-\frac{\pi}{2}}^{\frac{\pi}{2}}
\exp\Big\{-\frac{10}{3\pi^2}\,(\lambda-\lambda_2)\,\theta^2\Big\}\,\theta^6\,d\theta \\
 & \leq &
\frac{1}{\pi}\, \Big(\frac{20}{3\pi^2}\,(\lambda-\lambda_2)\Big)^{-7/2}
\int_{-\infty}^\infty
e^{-x^2/2}\,x^6\,dx \\ 
 & = &
\pi^{\frac{13}{2}}\, \Big(\frac{3}{20}\Big)^{7/2}\, 15 \sqrt{2} \ 
\frac{1}{(\lambda-\lambda_2)^{7/2}}  \, < \, 
\frac{48}{(\lambda-\lambda_2)^{7/2}}.
\ene

Now, the assumption (5.3) may be rewritten as
$$
{\rm Im}\,\Big(\sum_{l=1}^n 
\log(q_l+p_lr\,e^{i\theta})\Big)'\Big|_{\theta=0}=
\Big(\sum_{l=1}^n 
{\rm Im}\big(\log(q_l+p_lr\,e^{i\theta})\big)\Big)'\Big|_{\theta=0}=k.
$$
Here, the functions 
${\rm Im}\big(\log(q_l + p_l r\,e^{i\theta})\big)$ are odd, so their
2nd derivatives are vanishing at zero. We now apply 
the Taylor formula up to the cubic term to the function
$$
A_k(r,\theta) = -k\theta + {\rm Im}\, \sum_{l=1}^n \log(q_l+p_lre^{i\theta})
$$
on the interval $\theta\in [-\pi/2,\pi/2]$ to get that 
$$
A_k(r,\theta) \, = \, \frac{1}{6}\, \Big({\rm Im}\, 
\sum_{l=1}^n \log(q_l+p_l r e^{iv})\Big)^{'''}\Big|_{v=\theta_0}\,\theta^3
$$
with some $\theta_0 \in [-\frac{\pi}{2},\frac{\pi}{2}]$. 
To perform differentiation, consider a function of the form
$$
h(v) = \log(q + p r\, e^{iv}), \qquad p,q,r>0.
$$
We have
\bee
h'(v)
 & = &
\frac{pr\, i e^{iv}}{q + p r\, e^{iv}} \, = \,
i\, \Big(1 - \frac{q}{q + p r\, e^{iv}}\Big) \, = \, 
i - iq\,(q + p r\, e^{iv})^{-1}, \\
h''(v) 
 & = &
-pqr\ e^{iv}\,(q + p r\, e^{iv})^{-2}, \\ 
h'''(v) 
 & = &
-pqr\ \Big(i e^{iv}\,(q + p r\, e^{iv})^{-2} - 
2i\,p r\, e^{2iv}\,(q + p r\, e^{iv})^{-3}\Big).
\ene
Therefore,
$$
-\Big({\rm Im}\, \sum_{l=1}^n \log(q_l + p_l r\,e^{i\theta})\Big)^{'''}
 \, = \,
{\rm Im}\, \Big(i \sum_{l=1}^n 
\frac{p_l q_l r\, e^{i\theta}}{(q_l + p_l r\, e^{i\theta})^2}\Big) - 2\,
{\rm Im}\, \Big(i \sum_{l=1}^n
\frac{q_l p_l^2\, r^2\,e^{2i\theta}}{(q_l + p_l r\, e^{i\theta})^3}\Big),
$$
implying that
$$
\bigg|\Big({\rm Im}\, \sum_{l=1}^n \log(q_l + p_l r\,e^{i\theta})\Big)^{'''}\bigg|
 \ \leq \
\sum_{l=1}^n \frac{p_l q_l r}{|q_l + p_l r\, e^{i\theta}|^2} + 2\,
\sum_{l=1}^n \frac{q_l p_l^2\, r^2}{|q_l + p_l r\, e^{i\theta}|^3}.
$$
But, for $\frac{5}{6} \le r\le 1$ and $|\theta|\le\frac{\pi}{2}$,
\bee
|q_l+p_lre^{i\theta}|^2 
 & = & 
(q_l+p_lr)^2 
(1-\frac{4q_lp_l r}{(q_l+p_lr)^2}\, \sin^2\frac{\theta}{2}\Big) \\
 & \geq &
(q_l+p_lr)^2 - 2q_lp_l r  \, = \, q_l^2 + p_l^2 r^2.
\ene
Hence
$$
\frac{r}{|q_l + p_l r\, e^{i\theta}|^2} \leq 
\frac{r}{q_l^2 + p_l^2 r^2} = u _l(r) \leq \frac{121}{60}.
$$

Here we used the property that $u_l(r)$ is increasing in 
$r \leq r_l = q_l/p_l$ and is decreasing in $r \geq r_l$.
If $r_l \geq 1$, this gives $u_l(r) \leq u_l(1) = \frac{1}{q_l^2 + p_l^2} \leq 2$.
If $r_l \leq \frac{5}{6}$, that is, when $p_l \geq \frac{6}{11}$, we get 
$u_l(r) \leq u_l(5/6) = \frac{5/6}{q_l^2 + \frac{5}{6}\, p_l^2}$.
The latter expression is minimized at $p_l = \frac{6}{11}$
where it has the value $\frac{121}{66}$.
Finally, if $\frac{5}{6} \leq r_l \leq 1$, which is equivalent to 
$ \frac{1}{2} \leq p_l \leq \frac{6}{11}$, we have
$$
u_l(r) \leq u_l(r_l) = \frac{1}{2p_l q_l} \leq 
\frac{1}{2 \cdot \frac{6}{11} \cdot \frac{5}{11}} = \frac{121}{60}.
$$

From this,
$$
\frac{r^2}{|q_l + p_l r\, e^{i\theta}|^3} \leq 
\Big(\frac{r^{4/3}}{q_l^2 + p_l^2 r^2}\Big)^{3/2} \leq
\Big(\frac{r}{q_l^2 + p_l^2 r^2}\Big)^{3/2} = 
u_l(r)^{3/2} \leq \Big(\frac{121}{60}\Big)^{3/2},
$$
so that
$$
\bigg|\Big({\rm Im}\, 
\sum_{l=1}^n \log(q_l + p_l r\,e^{i\theta})\Big)^{'''}\bigg|
 \, \leq \,
\frac{121}{60}
\sum_{l=1}^n p_l q_l + 2\,\Big(\frac{121}{60}\Big)^{3/2}\,
\sum_{l=1}^n q_l p_l^2  \, \leq \, c_0\,(\lambda-\lambda_2)
$$
with $c_0 = \frac{121}{60}+ 2\,(\frac{121}{60})^{3/2} < 7.744438$.
Thus, 
$$
|A_k(r,\theta)| \, \leq \, \frac{c_0}{6}\,(\lambda-\lambda_2)\,|\theta|^3, 
\qquad
\frac{5}{6} \leq r \leq 1, \ |\theta| \leq \frac{\pi}{2}.
$$

Now, as we mentioned before, the function $A_k$ is odd in $\theta$, so that
$I_{k1}(r)$ is a real number given by
\bee
I_{k1}(r)
 & = &
\frac{1}{2\pi} \int_{-\frac{\pi}{2}}^{\frac{\pi}{2}}\, \prod_{l=1}^n
\frac{|q_l + p_l r\, e^{i\theta}|}{q_l + p_l r}\,
\cos(A_k(r,\theta))\,d\theta \\
 & = &
\frac{1}{2\pi} \int_{-\frac{\pi}{2}}^{\frac{\pi}{2}}\, \prod_{l=1}^n
\frac{|q_l+p_lre^{i\theta}|}{q_l + p_l r}\,d\theta
-\frac{1}{\pi} \int_{-\frac{\pi}{2}}^{\frac{\pi}{2}}\, \prod_{l=1}^n
\frac{|q_l+p_lre^{i\theta}|}{q_l + p_lr}\,
\sin^2(A(r,\theta)/2)\,d\theta.
\ene
Hence, using
$$
\sin^2(A(r,\theta)/2) \, \leq \, \frac{1}{4} A_k(r,\theta)^2 \, \leq \, 
\frac{c_0^2}{144}\, (\lambda-\lambda_2)^2\,\theta^6,
$$
from the previous estimates we may deduce the lower bound
\bee
I_{k1}(r)
 & \ge & 
0.3093\,\frac{1}{\sqrt{\lambda-\lambda_2}} - 
\frac{c_0^2}{144}\,(\lambda-\lambda_2)^2 \, 
\frac{48}{(\lambda-\lambda_2)^{7/2}} \\
 & = & 
0.3093\,\frac{1}{\sqrt{\lambda-\lambda_2}} - \frac{c_0^2}{3}\,
\frac{1}{(\lambda-\lambda_2)^{3/2}} \\
 & \geq &
\frac{1}{\sqrt{\lambda-\lambda_2}} \, 
\Big(0.3093 - \frac{20}{\lambda-\lambda_2}\Big) \, \geq \, 
0.1093\, \frac{1}{\sqrt{\lambda-\lambda_2}},
\ene
where on the last step we assume that $\lambda - \lambda_2 \geq 100$. 
Together with the upper bound on $I_{k2}$, we arrive at the lower bound
\bee
I_k(r) 
 & \ge &
0.1093\, \frac{1}{\sqrt{\lambda-\lambda_2}} - \frac{1}{2}\,
e^{-\frac{5}{6}\,(\lambda-\lambda_2)} \\ 
 & \ge &
\big(0.1093 - 5\, e^{-\frac{500}{6}}\big)\,
\frac{1}{\sqrt{\lambda-\lambda_2}}
 \, > \,\frac{0.1}{\sqrt{\lambda-\lambda_2}}.
\ene
Thus, Lemma 5.4 is proved. 
\qed

\vskip5mm
{\bf Proof of Proposition 5.1.}
We conclude from Lemmas 5.3 and 5.4 that
\be
\P\{W=k\} \, \ge \, \frac{1}{10\sqrt{\lambda-\lambda_2}}\,
e^{-4\,\frac{(\lambda-k)^2}{\lambda-\lambda_2}}
\en
for $0 \le \lambda - k \le \frac{1}{6}\,(\lambda-\lambda_2)$ under the assumption 
$\lambda-\lambda_2\ge 100$.

On the other hand, $f(k) = \P\{Z=k\} \le \frac{1}{\sqrt{2\pi k}}$, 
cf. (3.5). Since 
$k \geq \lambda - \frac{1}{6}\,(\lambda-\lambda_2) \geq \frac{5}{6}\,\lambda$,
we have 
$$
f(k) \, \le \, \frac{\sqrt{6/5}}{\sqrt{2\pi \lambda}} < 
\frac{1}{2\sqrt{\lambda}}.
$$
As a consequence,
\bee
1 + \chi^2(W,Z)
 & \ge &
\sum_{0 \le \lambda - k \le \frac{1}{6}\sqrt{\lambda-\lambda_2}}
\frac{\P\{W=k\}^2}{f(k)} \\
 & \geq &
\frac{\sqrt{\lambda}}{50\,(\lambda-\lambda_2)} \
\sum_{0 \le \lambda - k \le \frac{1}{6} \sqrt{\lambda-\lambda_2}} \,
e^{-8\,\frac{(\lambda-k)^2}{\lambda-\lambda_2}} \ \ge \
0.001\, \sqrt{\frac{\lambda}{\lambda-\lambda_2}}.
\ene
In order to clarify the last inequality, 
note that the condition $\lambda - \lambda_2 \geq 100$ implies that 
$\lambda > 100$. The above summation is performed over all integers 
$k$ from the interval 
$\lambda - \frac{1}{6}\sqrt{\lambda-\lambda_2} \leq x \leq \lambda$ of 
length at least $10/6$. It contains at least one integer point, and actually,
the number of integer points in it is at least
$h = \frac{1}{6}\sqrt{\lambda-\lambda_2}$. Moreover,
\bee
\sum_{0 \le \lambda - k \le h} \,
e^{-8\,\frac{(\lambda-k)^2}{\lambda-\lambda_2}} 
 & \geq & 
\sum_{[\lambda - h] + 1 \leq k \leq [\lambda]} \,
\int_{\lambda - k}^{\lambda - k + 1}
e^{-\frac{8\,x^2}{\lambda-\lambda_2}}\, dx  \\
 & = &
\int_{\lambda - [\lambda]}^{\lambda - [\lambda - h]} 
e^{-\frac{8 x^2}{\lambda-\lambda_2}}\,dx \, \geq \,
\frac{1}{4}\sqrt{\lambda-\lambda_2} \int_{2/5}^{2/3} e^{-y^2/2}\,dy \\
 & = & 
\frac{\sqrt{2\pi}}{4}\sqrt{\lambda-\lambda_2}\ \big(\Phi(2/3) - \Phi(2/5)\big)
 \, \geq \, 0.056\sqrt{\lambda-\lambda_2}.
\ene
Here, we used the bounds 
$4\,\frac{\lambda - [\lambda]}{\sqrt{\lambda-\lambda_2}} \leq \frac{2}{5}$
and 
$4\,\frac{\lambda - [\lambda - h]}{\sqrt{\lambda-\lambda_2}} \geq 
4\,\frac{\lambda - [\lambda - 10/6]}{10} \geq \frac{2}{3}$, 
together with $\Phi(2/3) - \Phi(2/5) > 0.09$.

In order to treat the region $\lambda-\lambda_2\le 100$, 
we apply Proposition 2.2. Let $W_1 = W$ and $W_2=Y_1+\dots+Y_{m}$, 
where $Y_1,\dots Y_m$ are independent Bernoulli random variables taking 
values 1 and 0 with probabilities $1/2$ and $m = 400$. Assume as well 
that $W$ and $W_2$ are independent. Then
$\tilde{\lambda} = \lambda+m/2$ and $\tilde{\lambda}_2 = \lambda_2 + m/4$ 
satisfy the condition $\tilde{\lambda}-\tilde{\lambda}_2 \ge 100$.

Denote by $Z_2$ a Poisson random variable with $\E Z_2 = m/2$ which is 
independent of $Z_1 = Z$. By the previous step and the inequality (2.4)
of Proposition 2.2,
\bee
0.001\,\sqrt{\frac{\tilde{\lambda}}{\tilde{\lambda}-\tilde{\lambda}_2}} 
 & \le &
\chi^2(W_1+W_2,Z_1+Z_2) + 1 \\
 & \le &
(\chi^2(W_1,Z_1)+1)(\chi^2(W_2,Z_2)+1).
\ene
Here, by (4.1), $\chi^2(W_2,Z_2) \leq 19\sqrt{2}$. Moreover, since
$\lambda-\lambda_2\le 100$, we have
$$
\sqrt{\frac{\tilde{\lambda}}{\tilde{\lambda}-\tilde{\lambda}_2}} \, = \,
\sqrt{\frac{\lambda + m/2}{\lambda - \lambda_2 + m/4}} \, \geq \,  
\sqrt{\frac{\lambda + 200}{200}} \, \geq \, \frac{1}{10\sqrt{2}}\,
\sqrt{\frac{\lambda}{\max\{1,\lambda - \lambda_2\}}}\, .
$$
It follows that
$$
1+\chi^2(W,Z) \geq \frac{0.001}{10\sqrt{2}\,(19\sqrt{2} + 1)}
\sqrt{\frac{\lambda}{\max\{1,\lambda - \lambda_2\}}} > 2.5 \cdot 10^{-6}\,
\sqrt{\frac{\lambda}{\max\{1,\lambda - \lambda_2\}}}\, .
$$
Hence, Proposition 5.1 holds in the case $\lambda-\lambda_2\le 100$ as well.

%---------------------------- Section 6 ----------------------------------------
\vskip7mm
\section{{\bf \large Lower Bound on $D$}}
\setcounter{equation}{0}

\vskip2mm
\noindent
An analogue of Proposition 5.1 is the following statement for the relative
entropy. Recall that $Q = \lambda/\max\{1,\lambda-\lambda_2\}$.

\vskip5mm
{\bf Proposition 6.1}. {\sl 
If $\lambda_2 \geq \kappa_0 \lambda$ and $\lambda\ge \lambda_0$, then
\be
D(W||Z) \, \ge \, c_0 \log(eQ),
\en
where $\kappa_0 = 1 - \exp\{-2\cdot 10^7\}$, 
$\lambda_0 = \exp\{2\cdot 10^7\}$, and $c_0 = e^{-14}$.
}

\vskip5mm
{\bf Proof.} 
Let us recall two estimates from the previous section, namely
\bee
w_k 
 & = &
\P\{W=k\} \, \ge \, \frac{1}{10\sqrt{\lambda-\lambda_2}}\,
e^{-4\,\frac{(\lambda-k)^2}{\lambda-\lambda_2}}, \\
v_k 
 & = & 
\P\{Z=k\} \, \le \, \frac{1}{\sqrt{2\pi k}}.
\ene
The first one is valid under the conditions
$0 \le \lambda - k \le \frac{1}{6}\,(\lambda-\lambda_2)$ and
$\lambda-\lambda_2\ge 100$, cf. (5.4). Clearly, they are fulfilled if 
$0 \le \lambda - k \le \frac{5}{3}\sqrt{\lambda-\lambda_2}$ and 
$\lambda-\lambda_2\ge 100$. If additionally $\lambda_2 \geq \kappa \lambda$, 
$0 < \kappa < 1$, then
$$
w_k \, \ge \, \frac{1}{10\sqrt{\lambda-\lambda_2}}\, e^{-100/9} \, \ge \, 
\frac{1}{10\sqrt{(1 - \kappa)\,\lambda}}\, e^{-100/9}.
$$
Since $k \geq \frac{5}{6}\,\lambda$, we also have an upper bound
$$
v_k \, \le \, \frac{1}{\sqrt{5\pi \lambda/3}}.
$$
In order that $w_k \geq v_k$, it is therefore sufficient to require that
$\frac{1}{10\sqrt{1 - \kappa}}\, e^{-100/9} \geq \frac{1}{\sqrt{5\pi/3}}$,
that is, $1 - \kappa \leq \frac{\pi}{60}\,e^{-200/9}$.
We have, moreover,
$$
\log \frac{w_k}{v_k}  \, \ge \, 
\frac{1}{2}\log\frac{e\lambda}{\lambda-\lambda_2} 
+ \log\Big(\frac{\sqrt{5\pi/3e}}{10}\, e^{-100/9}\Big) \, \ge \, 
\frac{1}{2}\log\frac{e\lambda}{\lambda-\lambda_2} - 14.
$$
Now, applying the inequality (2.1) of Proposition 2.1, we get
\bee
D(W||Z)
 & \ge &
\sum_{w_k \geq v_k} w_k \log \frac{w_k}{v_k} - 1 \\
 & \ge &
\sum_{0 \leq \lambda - k \le \frac{5}{3}\sqrt{\lambda-\lambda_2}} 
w_k \log \frac{w_k}{v_k} - 1 \\
 & \ge &
\sum_{0 \le \lambda - k \le \frac{5}{3}\sqrt{\lambda-\lambda_2}} 
w_k\,\Big(\frac{1}{2}\log\frac{e\lambda}{\lambda-\lambda_2} - 14\Big) - 1 \\ 
 & \ge &
\frac{1}{2} \log\frac{e\lambda}{\lambda-\lambda_2}
\sum_{0 \le \lambda - k \le \frac{5}{3}\sqrt{\lambda-\lambda_2}} 
\frac{1}{10\sqrt{\lambda-\lambda_2}}\,
e^{-4\,\frac{(\lambda-k)^2}{\lambda-\lambda_2}} - 15.
\ene
Note that, if $\lambda - \lambda_2 \geq 100$, the $x$-interval 
$0 \le \lambda - x \le \frac{5}{3}\sqrt{\lambda-\lambda_2}$ 
has length at least $50/3$, so, the total number of integer points 
in this interval is at least $50/3$ as well.
Hence, the last sum can be bounded from below by
$$
\frac{50/3}{10\sqrt{\lambda-\lambda_2}}\, e^{-100/9}
\sum_{0 \le \lambda - k \le \frac{5}{3}\sqrt{\lambda-\lambda_2}} 1 \, \geq \,
\frac{5}{3}\,e^{-100/9} \, > \, e^{-11}.
$$
Thus,
\be
D(W||Z) \, \geq \,
\frac{1}{2}\, e^{-11} \log\frac{e\lambda}{\lambda-\lambda_2} - 15.
\en
Moreover, if $\lambda_2 \geq \kappa \lambda$ with 
$\kappa \geq \kappa_1 = 1-\exp\{- 60 \, e^{11}\}$, then
$$
\frac{1}{4}\, e^{-11}\, \log\frac{e\lambda}{\lambda-\lambda_2} \, \geq \, 
\frac{1}{4}\, e^{-11}\, \log\frac{1}{1-\kappa}
\, \geq \, 15,
$$
and (6.2) yields
\be
D(W||Z) \, \geq \,
\frac{1}{4}\, e^{-11}\, \log\frac{e\lambda}{\lambda-\lambda_2}.
\en

The proposition is thus proved under the conditions 
$\lambda-\lambda_2\ge 100$ and $\lambda_2 \geq \kappa \lambda$ with
$\kappa_1 \leq \kappa < 1$.
It remains to eliminate the first condition, assuming that 
$\lambda-\lambda_2 < 100$ and again that $\lambda_2\ge \kappa \lambda$
with $\kappa$ being sufficiently close to 1.
To this aim, we appeal to Proposition 2.2 again like in the last step of the
proof of Proposition 5.1. Namely, using the same notations
and assumptions, from the inequality (2.3) and using (6.3),
we obtain that
\begin{eqnarray}
\frac{1}{4}\, e^{-11}\,
\log\frac{e\tilde \lambda}{\max\big\{1,\tilde \lambda-\tilde \lambda_2\big\}} 
 & \le &
D(W_1+W_2||Z_1+Z_2) \nonumber \\
 & \le &
D(W_1||Z_1) + D(W_2||Z_2),
\end{eqnarray}
where $W_1 = W$ and $Z_1 = Z$. It holds, as long as
$\tilde{\lambda}_2 \ge \kappa \tilde{\lambda}$, i.e.,
$
\lambda_2+m/4 \, \geq \, \kappa\,(\lambda + m/2).
$
Since $\lambda-\lambda_2 < 100$, the latter would follow from
$$
\lambda - 100 + m/4 \, \geq \, \kappa\,(\lambda + m/2)
$$
which is solved as 
$$
\lambda \geq 200\,\frac{\kappa}{1 - \kappa}.
$$

Moreover, by (4.2), we have $D(W_2||Z_2) \leq 23\,\log(2e)$. This
bound may be used in (6.4), which gives
\bee
D(W||Z) 
 & \geq & 
\frac{1}{4}\, e^{-11}\,
\log\frac{e\tilde \lambda}{\max\big\{1,\tilde \lambda-\tilde \lambda_2\big\}} - 
23\,\log(2e) \\
 & \geq & 
\frac{1}{8}\, e^{-11}\,
\log\frac{e\tilde \lambda}{\max\big\{1,\tilde \lambda-\tilde \lambda_2\big\}},
\ene
where the second inequality holds true when $1-\kappa$ is sufficiently
small. Namely,
$$
\frac{1}{8}\, e^{-11}
\log\frac{e\tilde \lambda}{\tilde \lambda-\tilde \lambda_2} \, \geq \, 
\frac{1}{8}\, e^{-11} \log\frac{1}{1-\kappa}
\, \geq \, 23\,\log(2e),
$$
if $\tilde \lambda_2 \geq \kappa \tilde \lambda$ and
$1 - \kappa \leq \exp\{- 8 \cdot 23 \cdot \log(2e) \cdot e^{11}\}$.
Since the product in the exponent is smaller than $1.87 \cdot 10^7$,
we may choose $\kappa = 1-\exp\{- 1.87 \cdot 10^7\} > \kappa_1$. In this case,
$$
D(W||Z) \, \geq \, c_1\, 
\log\frac{e\tilde \lambda}{\tilde \lambda - \tilde \lambda_2}, \qquad
c_1 = \frac{1}{8}\, e^{-11},
$$
assuming that $\lambda \geq 200\,\frac{\kappa}{1 - \kappa}$. 
But
$$
\log\frac{e\tilde \lambda}{\tilde \lambda - \tilde \lambda_2} =
\log\frac{e\,(\lambda + 200)}{\lambda - \lambda_2 + 100} \geq
\frac{1}{2}\, \log\frac{e\lambda}{\max\{1,\lambda - \lambda_2\}}
$$
for all $\lambda \geq 4 \cdot 10^4$.
It remains to note that $200\,\frac{\kappa}{1 - \kappa} < \lambda_0$, 
$\kappa < \kappa_0$, $\frac{1}{2}\, c_1 > c_0$.
\qed

%---------------------------- Section 7 ----------------------------
\vskip7mm
\section{{\bf \large Proof of Theorem 1.1}}
\setcounter{equation}{0}

\vskip2mm
\noindent
Let us summarize. Using the quantity
$$
F = F(\lambda,\lambda_2) = 
\frac{\max(1,\lambda)}{\max(1,\lambda-\lambda_2)},
$$
the results on Poisson approximation obtained 
for different regions of $\lambda$ and $\lambda_2$ can be combined
in the form of the following two-sided bounds
\be
c_1\, \Big(\frac{\lambda_2}{\lambda}\Big)^2\, (1 + \log F)
 \, \le \, D(W||Z) \, \le \,
c_2\Big(\frac{\lambda_2}{\lambda}\Big)^2\, (1 + \log F),
\en
\be
c_1\,\Big(\frac{\lambda_2}{\lambda}\Big)^2 \sqrt{F} \, \le \, \chi^2(W,Z)
 \, \le \,
c_2\,\Big(\frac{\lambda_2}{\lambda}\Big)^2 \sqrt{F},
\en
which are valid up to some absolute positive constants $c_1$ and $c_2$.
Let us describe the proof of Theorem 1.1 and provide explicit values for
these constants. As we will see, (7.1)-(7.2) hold 
with $c_1 = 10^{-8}$ and $c_2 = 5.6 \cdot 10^7$.

\vskip2mm
{\bf An upper bound in} (7.1).

If $\lambda\le 1/2$, these bounds simplify and are made precise via
\be
\frac{1}{4}\,\Big(\frac{\lambda_2}{\lambda}\Big)^2 \, \le \, 
D(W||Z) \, \leq \, \chi^2(W,Z) \, \le \,
15\,\Big(\frac{\lambda_2}{\lambda}\Big)^2.
\en
Here, the left inequality holds for all $\lambda$ and $\lambda_2$, 
cf. [H-J-K], while the right inequality
is part of Proposition 3.1. Note that 
$\lambda\le 1/2$ implies $\lambda_2 \le \frac{1}{2}\, \lambda$.

If $\lambda \geq 1/2$ and $\lambda_2 \le \frac{1}{2}\, \lambda$,
we have, by Proposition 3.2, 
$$
D(W||Z) \, \leq \, \chi^2(W,Z) \, \le \,
56 \cdot 10^6\,\Big(\frac{\lambda_2}{\lambda}\Big)^2,
$$
so that
\be
D(W||Z) \, \leq \,
56 \cdot 10^6\,\Big(\frac{\lambda_2}{\lambda}\Big)^2\,(1 + \log F).
\en

In the case where $\lambda \geq 1/2$ and $\lambda_2 > \frac{1}{2}\, \lambda$,
one may apply (4.2) which gives
$$
D(W||Z) \, \leq \, 23\,(1 + \log F) \, \leq \,
4 \cdot 23\,\Big(\frac{\lambda_2}{\lambda}\Big)^2\,(1 + \log F).
$$
Here, the right-hand side contains a better numerical constant
in comparison with (7.4), and we finally get (7.1) with a constant
$c_2 = 56 \cdot 10^6$.

\vskip2mm
{\bf A lower bound in} (7.1).

If $\lambda\le 1$, then $F = 1$, so that the lower bound in (7.3) yields
(7.1) with $c_1 = 1/4$.

If $\lambda \geq 1$, the inequality (7.4) may be reversed by virtue of 
(6.1), which gives
\be
D(W||Z) \, \geq \, c_0 (1+\log F) \, \geq \, 
c_0 \Big(\frac{\lambda_2}{\lambda}\Big)^2\,(1+\log F)
\en
with $c_0 = e^{-14}$, provided that 
$\lambda_2 \geq \kappa_0 \lambda$ and $\lambda \geq \lambda_0$, where
$\kappa_0 = 1 - \exp\{-2\cdot 10^7\}$ and $\lambda_0 = \exp\{2\cdot 10^7\}$.
But, the remaining regions belong to the non-degenerate case, where 
$F$ is bounded by a quantity which depends on $\kappa_0$ or $\lambda_0$.
Indeed, if $\lambda_2 \leq \kappa_0 \lambda$, then 
$\log F \leq -\log(1 - \kappa_0) = 2\cdot 10^7$, so,
$$
D(W||Z) \, \geq \, \frac{1}{4\,(1 + 2\cdot 10^7)}\,
\Big(\frac{\lambda_2}{\lambda}\Big)^2\,(1 + \log F).
$$
This means that the left inequality in (7.1) holds with a constant
$c_1 = \frac{1}{4\,(1 + 2\cdot 10^7)}$ which is smaller than $c_0$
in the analogous inequality (7.5). Similarly, if
$1 \leq \lambda < \lambda_0$, then
$F \leq \lambda < \lambda_0$, and we get, by the lower bound in (7.3),
$$
D(W||Z) \, \geq \, \frac{1}{4\,(1 + \log \lambda_0)}\,
\Big(\frac{\lambda_2}{\lambda}\Big)^2\,(1 + \log F).
$$
This means that the left inequality in (7.1) holds true with the same 
constant $c_1$ as above. Thus, the lower bound in (7.1) holds with 
constant $c_1$ ($> 10^{-8}$).

\vskip2mm
{\bf An upper bound in} (7.2).

If $\lambda\le 1/2$, we have (7.3), which implies (7.2) with $c_2 = 15$.
 
If $\lambda \geq 1/2$ and $\lambda_2 \le \frac{1}{2}\, \lambda$, 
a stronger version of (7.4) is provided by Proposition 3.2, which gives
$$
\chi^2(W,Z) \, \le \,
56 \cdot 10^6\,\Big(\frac{\lambda_2}{\lambda}\Big)^2,
$$
so that (7.2) holds true with $c_2 = 56 \cdot 10^6$.
In the case where $\lambda \geq 1/2$ and $\lambda_2 > \frac{1}{2}\, \lambda$,
one may apply (4.1) which gives
$$
\chi^2(W,Z) \, \leq \,
76\,\Big(\frac{\lambda_2}{\lambda}\Big)^2\sqrt{F}.
$$
Here, the right-hand side contains a better numerical constant, and 
we finally get (7.2) with the same constant $c_2$ as in (7.1).

\vskip2mm
{\bf A lower bound in} (7.2).

If $\lambda\le 1$, then $F = 1$, so that the lower bound in (7.3) yields
(7.1) with $c_1 = 1/4$.

Assume that $\lambda \geq 1$, in which case
$F = Q = \lambda/\max(1,\lambda - \lambda_2)$. By (5.2), we have
$$
\chi^2(W,Z) \, \geq \, \frac{c_0}{9}\sqrt{F}
$$
with $c_0 = 2.5 \cdot 10^{-6}$, provided that 
$\lambda_2 \geq \kappa_0 \lambda$, $\kappa_0 = 1 - c_0^2/4$.
This gives
\be
\chi^2(W,Z) \, \geq \, \frac{c_0}{9}
\Big(\frac{\lambda_2}{\lambda}\Big)^2 \sqrt{F},
\en
and we obtain the left inequality in (7.2) with $c_1 = c_0/9 > 10^{-7}$.

The remaining region belongs to the non-degenerate case, where $F$ is 
bounded. Indeed, if $\lambda_2 \leq \kappa_0 \lambda$, then 
$1/\sqrt{F} \geq \sqrt{1 - \kappa_0} = \frac{c_0}{2} = 0.8\,\cdot 10^{-6}$, 
so that, by the left inequality in (7.3),
$$
\chi^2(W,Z) \, \geq \, \frac{1}{4}\,
\Big(\frac{\lambda_2}{\lambda}\Big)^2 \, \geq \, 
0.2\,\cdot 10^{-6}\,\Big(\frac{\lambda_2}{\lambda}\Big)^2 \sqrt{F}.
$$
This means that the left inequality in (7.1) holds true with constant
$c_1 = 2\,\cdot 10^{-7}$ which is slightly better than the constant
in the analogous inequality (7.6). 
Thus, the lower bound in (7.2) holds true with constant $c_1 = 10^{-7}$.
\qed

%---------------------------- Section 8 ----------------------------
\vskip7mm
\section{{\bf\large Tsallis versus Vajda-Pearson}} 
\setcounter{equation}{0}

\vskip2mm
\noindent
We now turn to the Tsallis relative entropies of other indexes. 
To make an application of non-uniform bounds more convenient, 
first let us relate $T_\alpha$ to the Vajda-Pearson distance
$$
\chi_\alpha(X,Z) = \int \Big|\frac{p-q}{q}\Big|^\alpha\,q\,d\pi.
$$
It is defined for arbitrary random elements $X$ and $Z$ in 
a measure space $(\Omega,\pi)$ whose distributions are
absolutely continuous and have densities $p$ and $q$
respectively with respect to the measure $\pi$ on $\Omega$
(the defnition does not depend on the choice of $\pi$).

Recall that
$$
T_\alpha(X||Z) = \frac{1}{\alpha - 1}\,\bigg[
\int \Big(\frac{p}{q}\Big)^\alpha\,q\,d\pi - 1\bigg],
$$
so that $T_2 = \chi_2$ is the classical Pearson distance, and
note that $T_\alpha = \chi_\alpha = \infty$ as long as
the distribution of $X$ is not absolutely continuous with
respect to the distribution of $Z$.
We need the following auxilliary result.

\vskip5mm
{\bf Proposition 8.1.} {\sl For any $\alpha\ge 2$,
$$
T_{\alpha}(W||Z) \, \le \, \frac{2^{\alpha}}{\alpha-1}\,
\Big(T_{2}(W||Z) + \chi_{\alpha}(W,Z)\Big).
$$
}

{\bf Proof.} We may assume that the distribution of $X$ 
is absolutely continuous with respect to the distribution 
of $Z$, with $\chi_{\alpha}(W,Z) < \infty$. In this case, 
the (non-negative) function $\xi = p/q$ is well defined a.e. 
with respect to the probability measure $Q = q\,d\pi$.
We consider it as a random variable on the probability space
$(\Omega,Q)$ with finite moment of order $\alpha$.
Note that 
$$
(\alpha-1)\,T_{\alpha}(W||Z) = 
\E\, (\xi^{\alpha}-1) \quad {\rm and} \quad 
\chi_{\alpha}(W,Z) = \E\,|\xi-1|^{\alpha}. 
$$
Putting $\eta = \xi-1 \ge -1$, define the function 
$\psi(t) = \E\, (1+t\eta)^{\alpha} -1$, $t \ge 0$, so that
$\psi(1) = (\alpha-1)\,T_{\alpha}(W||Z)$. 
By the integral Taylor formula,
$$
\psi(1) \, = \, \alpha(\alpha-1)\, \E\, \eta^2 \int_{0}^{1}
(1-t)(1+t\eta)^{\alpha-2}\,dt.
$$
Introducing the sets $A = \{\xi\le 2\} = \{\eta \le 1\}$ and 
$B = \{\xi>2\} = \{\eta > 1\}$, we have
\bee
\E \,1_A\,\eta^2\int_0^1 (1-t)(1+t\eta)^{\alpha-2}\,dt 
 & \le &
\E\, \eta^2\int_0^1 (1-t)(1+t)^{\alpha-2}\,dt \\
 & \le & 
\frac{2^{\alpha}}{\alpha(\alpha-1)}\ T_2(W||Z)
\ene
and
\bee
\E\, 1_B \, \eta^2\int_{0}^{1}(1-t)(1+t\eta)^{\alpha-2}\,dt
 & \le & 
\E\, 1_B \, \eta^{\alpha}\int_{0}^{1}(1-t)(1+t)^{\alpha-2}\,dt \\
 & \le &  
\frac{2^{\alpha}}{\alpha(\alpha-1)}\, \chi_{\alpha}(W,Z).
\ene
We obtain the assertion of the proposition from the last two bounds.
\qed

%---------------------------- Section 9 ----------------------------
\vskip7mm
\section{{\bf\large Estimates of Vajda-Pearson distances}}
\setcounter{equation}{0}

\vskip2mm
\noindent
For the proof of Theorem 1.2 we need the following propositions.
We thus return to the setting of Bernoulli trials. Let us denote by
$c(\alpha)$ a positive constant depending on $\alpha$ only,
which may vary from place to place.

\vskip5mm
{\bf Proposition 9.1.} {\sl For $\alpha > 1$ and $\lambda\le\frac{1}{2}$, 
we have
$$
\chi_{\alpha}(W,Z)\,\le\, 
c(\alpha)\,\frac{\lambda_2^{\alpha}}{\lambda^{2(\alpha-1)}}.
$$
}

%\vskip5mm
{\bf Proof.} Applying Lemmas III.1-2 and repeating the argument 
used in the proof of Proposition III.4 from \cite{B-C-G2}, we 
obtain that
\bee
\frac{e^{\lambda}}{\lambda_2^{\alpha}}\,\chi_{\alpha}(W,Z) 
 & \le & 
1 + \lambda\Big(\frac{\lambda+e-1}{\lambda}\Big)^{\alpha} \\
 &+ & 
3^{\alpha-1}\sum_{k=2}^{\infty}\,\frac{\lambda^k}{k!}\,
\Big(1+\Big(\frac{e^{\lambda}-1}{\lambda}\Big)^{\alpha}
\,\Big(\Big(\frac{k}{\lambda}\Big)^{\alpha}+
\Big(\frac{k(k-1)}{\lambda^2}\Big)^{\alpha}\Big)\\ 
&\le &
c(\alpha)\Big(1+\frac{1}{\lambda^{\alpha-1}} +
\frac{1}{\lambda^{2(\alpha-1)}}\Big).
\ene
\qed

\vskip5mm
{\bf Proposition 9.2.} {\sl Let $\alpha > 1$. If 
$\lambda\ge\, \frac{1}{2}$ and $\lambda_2\,\le\,\kappa\lambda$ with 
$\kappa\in(0,1)$, then
$$
\chi_{\alpha}(W,Z)\,\le \, \frac{c(\alpha)}{(1-\kappa)^{\frac{3\alpha}{2}}}\,
\Big(\frac{\lambda_2}{\lambda}\Big)^{\alpha}.
$$
}

{\bf Proof.} Write
$$
\chi_{\alpha}(W,Z) = \sum_{k=0}^{\infty}\,
\frac{|\Delta_k|^{\alpha}}{f(k)^{\alpha-1}}=S_1+S_2=
\Big(\sum_{k=0}^{[2\lambda]}+
\sum_{k=[2\lambda]+1}^{\infty}\,\Big)\,
\frac{|\Delta_k|^{\alpha}}{f(k)^{\alpha-1}}.
$$
In the range $0\le\, k\,\le [2\lambda]$ we apply the inequality 
(VI.2) from \cite{B-C-G2} which gives
$$
|\Delta_k|^{\alpha}\,\le\,\frac{c(\alpha)}{(1-\kappa)^{3\alpha/2}}\,
\Big(\frac{|k-\lambda|^{2\alpha}}{\lambda^{\alpha}}\,+\,1\Big)\,
\Big(\frac{\lambda_2}{\lambda}\Big)^{\alpha}\,f(k)^{\alpha}.
$$
Therefore
$$
S_1\le \,\frac{c(\alpha)}{(1-\kappa)^{\frac{3}{2}}}\,
\Big(\frac{\E\,|Z-\lambda|^{2\alpha}}{\lambda^{\alpha}}+1\Big)\,
\Big(\frac{\lambda_2}{\lambda}\Big)^{\alpha}\le\,
\frac{c(\alpha)}{(1-\kappa)^{\frac{3\alpha}{2}}}\,
\Big(\frac{\lambda_2}{\lambda}\Big)^{\alpha}.
$$
Here we use the upper bound 
$\E\,|Z-\lambda|^{2\alpha}\le c(\alpha)\,\lambda^{\alpha}$.

In order to estimate $S_2$ we use the inequalities (VI.3) and (II.1)
from \cite{B-C-G2} to get
\bee
S_2&\le & \frac{c(\alpha)}{(1-\kappa)^{\frac{3\alpha}{2}}}\,
\sum_{k=[2\lambda]+1}^{\infty}
\Big(\frac{k}{\lambda}\Big)^{3\alpha}\,\lambda_2^{\alpha}\,f(k)\\
&\le &
\frac{c(\alpha)}{(1-\kappa)^{\frac{3\alpha}{2}}}\,
\lambda_2^{\alpha} f([2\lambda]+1) \sum_{k=0}^{\infty}
\Big(1+\Big(\frac{k}{\lambda}\Big)^{3\alpha}\Big)\,
\frac{1}{2^{k}}\\
&\le &
\frac{c(\alpha)}{(1-\kappa)^{\frac{2\alpha}{2}}}\,
\lambda_2^{\alpha}\,f([2\lambda]+1)\,\le 
\frac{c(\alpha)}{(1-\kappa)^{\frac{3\alpha}{2}}}\,
\lambda_2^{\alpha}\,e^{\frac{1}{2}\lambda\log\frac{e}{4}}.
\ene
The assertion of the proposition follows immediately from 
the last two estimates.
\qed

%---------------------------- Section 10 ----------------------------
\vskip7mm
\section{{\bf \large Proof of Theorem 1.2}}
\setcounter{equation}{0}

\vskip2mm
\noindent
To complete the proof of Theorem 1.2, we need the following 
two lemmas. Recall that $Q=\lambda/\max\{1,\lambda-\lambda_2\}$. 

\vskip5mm
{\bf Lemma 10.1.} {\sl For $\alpha\,>1$ and $\lambda\ge\,\frac{1}{2}$,
$$
T_{\alpha}(W||Z) \, \le\, c(\alpha)\,Q^{(\alpha-1)/2}.
$$
}

\vskip2mm
{\bf Proof.}
By the definition of the Tsallis distance, 
\bee
T_{\alpha}(W||Z)
&\le &
\frac{1}{\alpha-1}\sum_{k=0}^{\infty}\,
\Big(\frac{w_k}{v_k}\Big)^{\alpha}\,v_k\ \le \,
S_1\,+\,S_2\,+\,S_3\\
& = &
\,\frac{1}{\alpha-1}\,\Big(\sum_{1\le k<\frac{1}{4}\lambda}
 \,+\, \sum_{\frac{1}{4}\lambda\le k\le 4\lambda}\,+\,
\sum_{k>4\lambda}\Big)\,
\Big(\frac{w_k}{v_k}\Big)^{\alpha-1}\,w_k.
\ene
By (4.5),
$$
(\alpha-1)S_2\,\le\,\sum_{\frac{1}{4}\lambda\le k\le 4\lambda}
\big(2e\,Q_0^{1/2}\,k^{1/2}\big)^{\alpha-1}\,w_k\,\le\,
(4e)^{\alpha-1}\,Q^{(\alpha-1)/2}.
$$
Using (4.6) and repeating the argument of Section 4, 
we obtain the upper bounds $S_1+S_3\,\le\, c(\alpha)$. 
The three last estimates give the assertion of the proposition.
\qed

\vskip5mm
{\bf Lemma 10.2.} {\sl For $\alpha\,>1$ and 
$\lambda\ge\,\frac{1}{2}$, with some constant $c_1(\alpha) \in (0,1)$
\be
1+T_{\alpha}(W||Z)\,\ge\,c(\alpha)\,Q^{(\alpha-1)/2}.
\en
Moreover
\be
T_{\alpha}(W||Z)\,\ge\frac{c_1(\alpha)}{9}\,Q^{(\alpha-1)/2}
\en
as long as $\lambda_2\,\ge\,\big(1-\frac{c_1(\alpha)^2}{4}\big)\,\lambda$.
}

\vskip5mm
{\bf Proof.} The assertion (10.2) follows from the assertion (10.1) 
in the same way as (5.2) follows from (5.1). Therefore we omit the proof.

In order to prove (10.1) we use the lower bound (5.4). Repeating 
the argument of the proof of Proposition 5.1, we easily obtain 
the lower bound, under the assumption $\lambda-\lambda_2\,\ge\, 100$, 
\bee
1+T_{\alpha}(W||Z)
& \ge &
\sum_{0\le \lambda-k\le\frac{1}{6}\sqrt{\lambda-\lambda_2}}\,
\Big(\frac{w_k}{v_k}\Big)^{\alpha-1}\,w_k\\
& \ge &
\frac{1}{10}\cdot\Big(\frac{1}{5}\Big)^{\alpha-1}\,
\Big(\frac{\lambda}{\lambda-\lambda_2}\Big)^{(\alpha-1)/2}\,
\frac{1}{\sqrt{\lambda-\lambda_2}}\sum_{0\le\lambda-k\le
\frac{1}{6}\sqrt{\lambda-\lambda_2}}\,
e^{-4\alpha\frac{(\lambda-k)^2}{\lambda-\lambda_2}}\\
&\ge&
c(\alpha)\Big(\frac{\lambda}{\lambda-\lambda_2}\Big)^{(\alpha-1)/2}.
\ene
In order to treat the region $\lambda-\lambda_2\,\le\, 100$ we refer 
to Johnson \cite{J}, pp. 133--134,  and repeat the argument of the end 
of Section 5.
\qed

\vskip5mm
{\bf Proof of Theorem 1.2.}
Assuming that 
$\frac{\lambda_2}{\lambda}\le 1-\frac{1}{4}c_1(\alpha)^2$, we have
$F\,\sim\, 1$ with involved constants depending on $\alpha$, and then
we need to show that 
$T_{\alpha}(W||Z) \sim (\frac{\lambda_2} {\lambda})^2$. 

In the case $1\,<\alpha\,\le 2$, we have
$$
\frac{1}{4}\,\Big(\frac{\lambda_2}{\lambda}\Big)^2\,\le 
T_{\alpha}(W||Z)\,\le T_2(W||Z)=\chi^2(W,Z).
$$ 

Turning to the case $\alpha\,\ge 2$, first let $\lambda\,\le\frac{1}{2}$. 
Since $\lambda_2\le \lambda^2$, by Propositions 8.1 and 9.1,
\bee
T_{\alpha}(W||Z)
 & \le &
c(\alpha)\big(T_2(W||Z) + \chi_{\alpha}(W,Z)\big)\\
 & \le & 
c(\alpha)\Big(T_2(W||Z)+
\frac{\lambda_2^{\alpha}}{\lambda^{2(\alpha-1)}}\Big)\,\le\, 
c(\alpha)\Big(\frac{\lambda_2}{\lambda}\Big)^2. 
\ene
Now, let $\lambda\,\ge\frac{1}{2}$. Then, by Propositions 8.1 and 9.2, 
we conclude that
\bee
T_{\alpha}(W||Z)&\le& c(\alpha)\big(T_2(W||Z) + 
\chi_{\alpha}(W,Z)\big) \\
 & \le & 
c(\alpha)
\Big(T_2(W||Z) + \Big(\frac{\lambda_2}{\lambda}\Big)^{\alpha}\Big)\,
\le\,c(\alpha)\Big(\frac{\lambda_2}{\lambda}\Big)^2.
\ene

It remains to consider the region
$\frac{\lambda_2}{\lambda}\,\ge \,1-\frac{1}{4}\,c_1(\alpha)^2$. 
But in this case, the assertion of the theorem immediately 
follows from Lemmas 10.1 and 10.2.

%----------------------------- Section 12 ----------------------------------
\vskip7mm
\section{{\bf \large Difference of Entropies}}
\setcounter{equation}{0}

\vskip2mm
\noindent
For the proof of Corollary 1.3, we shall use another functional
$$
H_2(Z) = \big(\E\,(\log v(Z))^2\big)^{1/2} = 
\Big(\sum_k v_k\, (\log v_k)^2\Big)^{1/2}, \qquad v_k = \P\{Z=k\},
$$
where $Z$ is an integer-valued random variable. Thus, while the Shannon
entropy $H(Z) = -\E\,\log v(Z)$ describes the average of the informational
content $-\log v(Z)$, the informational quantity $H_2(Z)$ represents 
the 2nd moment of this random variable.

An application of Theorem 1.1 is based upon the following elementary
relation. 

\vskip5mm
{\bf Proposition 11.1.} {\sl For all integer-valued random variables $W$ and 
$Z$ with finite entropies, we have
\be
H(W||Z) \, \leq \, \chi^2(W,Z) + H_2(Z)\sqrt{\chi^2(W,Z)}.
\en
}

%\vskip2mm
{\bf Proof.} We may assume that the distribution of $W$ is absolutely
continuous with respect to the distribution of $Z$ (since otherwise
$\chi^2(W,Z) = \infty$). Equivalently, for all $k \in \Z$,
$v_k = 0 \Rightarrow w_k = 0$, where $w_k = \P\{W=k\}$.
Define $t_k = w_k/v_k$ in case $v_k > 0$. Recalling the definition (1.8),
we then have
$$
H(W||Z) = 
\sum_{v_k>0} (t_k \log t_k)\,v_k + \sum_{v_k>0} (t_k - 1)\, v_k \log v_k.
$$
We now apply the inequality $t\log t \leq (t-1) + (t-1)^2$ ($t \geq 0$), 
obtaining
\bee
H(W||Z) 
 & \leq &
\sum_{v_k>0} (t_k-1)\,v_k + \sum_{v_k>0} (t_k-1)^2\,v_k + 
\sum_{v_k>0} (t_k - 1) v_k \log v_k \\
 & = &
\sum_k 
\frac{(w_k - v_k)^2}{v_k} + \sum_{v_k > 0} (w_k - v_k) \log v_k.
\ene
Here, the first sum in the last bound is exactly $\chi^2(W,Z)$, while, by 
Cauchy's inequality, the square of the last sum is bounded from above by
$$
\sum_k \frac{(w_k - v_k)^2}{v_k} \, \sum_k v_k\, (\log v_k)^2 \, = \,
\chi^2(W,Z)\,H_2^2(Z).
$$
\qed

%\vskip2mm
In view of (11.1), we also need:

\vskip5mm
{\bf Proposition 11.2.} {\sl If $Z$ has a Poisson distribution with parameter
$\lambda$, then
$$
H_2(Z) \leq 
\begin{cases}
\sqrt{50}\,\log(1+\lambda), \hskip5mm {\sl if} \ \  \lambda \geq 1, \\
5\sqrt{\lambda}\,\log(e/\lambda), \qquad {\sl if} \ \ \lambda \leq 1. \\
\end{cases}
$$
}

\vskip2mm
{\bf Proof.} Put $v_k = \P\{Z=k\}$. In particular,
$v_0 \, (\log v_0)^2 = \lambda^2\,e^{-\lambda}$ and
$v_1 \, (\log v_1)^2 = \lambda e^{-\lambda}\,(\lambda + \log(1/\lambda))^2$.
This shows that the above upper bound for small $\lambda$ 
can be reversed up to a constant. For $\lambda \leq 1$, given $k \geq 1$, from
$$
\log\frac{1}{v_k} = \lambda + \log k! + k\log\frac{1}{\lambda} 
\leq k^2\log\frac{e}{\lambda},
$$
we get
$$
\sum_{k \geq 1} v_k\, (\log v_k)^2 \leq \E\,Z^4
\log^2\Big(\frac{e}{\lambda}\Big) \leq 
24\,\lambda\,\log^2\Big(\frac{e}{\lambda}\Big).
$$
Hence, $H_2^2(Z) \leq 25\,\lambda\,\log^2(e/\lambda)$,
thus proving the second upper bound of the lemma.

Now, assuming that $\lambda \geq 1$, let us apply the lower bounds 
(3.6)-(3.7) from Lemma 3.3, which for all $k \geq 1$ give
$$
\log\frac{1}{v_k} \leq 1 + \frac{1}{2}\, \log k + \frac{1}{\lambda}\,(k - \lambda)^2
 \leq \log(ek) + \frac{1}{\lambda}\,(k - \lambda)^2
$$
and
$$
\log^2\frac{1}{v_k} \leq 2 \log^2(e(k+1)) + \frac{2}{\lambda^2}\,(k - \lambda)^4.
$$
Note that this bound is also true for $k=0$.
Using the concavity of the function $\log^2 x$ in $x \geq e$ and applying
Jensen's inequality, we therefore obtain that
\bee
\sum_{k = 0}^\infty v_k\, (\log v_k)^2 
 & \leq &
2\,\E\, \log^2(e(Z+1)) + \frac{2}{\lambda^2}\,\E\,(Z - \lambda)^4 \\
 & \leq &
2 \log^2(e(\lambda + 1)) + \frac{6\,(\lambda + 2)}{\lambda} \, \leq \, 
2\, \big(1 + \log(1+\lambda)\big)^2 + 18.
\ene
Hence $H_2(Z) \leq Cx$, $x = \log(1+\lambda) \geq \log 2$, with
$C^2 = 2\,(1 + \frac{1}{x})^2 +  \frac{18}{x^2} < 50$.

Applying the upper bound (3.6) from Lemma 3.3, we also see that this
upper bound on $H_2$ can be reversed up to a constant as well.

\qed

\vskip2mm
{\bf Remark 11.3.} With similar arguments, it follows that
$$
H(Z) \leq 
\begin{cases}
c\log(1+\lambda), \hskip7mm {\rm if} \ \  \lambda \geq 1, \\
c\lambda \log(e/\lambda), \qquad {\rm if} \ \ \lambda \leq 1, \\
\end{cases}
$$
which can be reversed modulo an absolute factor $c>0$. Hence,
$H_2(Z) \sim H(Z)$ as long as $\lambda$ stays bounded away from 
zero.

\vskip5mm
{\bf Proof of Corollary 1.3.} By Theorem 1.1 with $W$ as in (1.1) and
with a Poisson random variable $Z$ with parameter $\lambda$, we have
$$
\chi^2(W,Z) \, \leq \, 
C\,\Big(\frac{\lambda_2}{\lambda}\Big)^2\, \sqrt{2 + \lambda}
$$
up to some absolute constant $C$. Using this estimate in (11.1) and
applying Proposition 11.2, the desired inequality (1.9) immediately
follows (in view of $\lambda_2 \leq \lambda$).

To derive a more precise inequality illustrating the asymptotic behaviour
in $\lambda$ in the typical case $\lambda_2 \leq \frac{1}{2}\,\lambda$, 
let us apply once more Theorem 1.1 with its sharper bound
$$
\chi^2(W,Z) \, \leq \, 
C\,\Big(\frac{\lambda_2}{\lambda}\Big)^2,
$$
as in Proposition 3.1. By Proposition 11.1, this gives
$$
H(W||Z) \, \leq \, C\,(1 + H_2(Z))\,\frac{\lambda_2}{\lambda},
$$
It remains to note that $1 + H_2(Z) \leq C\log(2+\lambda)$.
according to Proposition 11.2.
\qed

\vskip5mm
{\bf Acknowledgement.} We would like to thank the referee for 
drawing our attention to the work by V. Zacharovas and H.-K. Hwang.
Thanks also to A. Zaitsev for drawing our attention to the work by 
I. S. Borisov and I. S.  Vorozheĭkin.

\end{document}